\documentclass[12pt]{amsart}
\usepackage{latexsym, amsmath, amscd, amssymb, amsthm} 
\newtheorem{te}{Theorem}

 \newtheorem{lm}{Lemma}
\begin{document}
 \title{On Differential equation   of  invariants of binary form} 
\author{L. Bedratyuk}\address{ Khmelnitskiy national university, Instytuts'ka st. 11, Khmelnitskiy , 29016, Ukraine}

\begin{abstract}
An explicit form of single first order PDE for invariants of  binary form are  found. By solving  the equation a minimal generation set for a ring of invariants and theirs syzygies are calculated in the cases $n\leq 6$ and $n=8.$ 
\end{abstract}
\email{bedratyuk@ief.tup.km.ua}
\maketitle

\section{Introduction}
\noindent
Let $V_n$ be a vector $k$-space of the binary forms of degree $n$ 
$$ 
u(x,y)=\sum_{i=0}^n \, \alpha_i { n \choose i } x^{n-i} y^i,
$$

\noindent
where  $\alpha_i \in k$, and  $k \,$  is a field of characteristic zero.
\noindent
Let us indentify the coordinate ring $R_n$  of the space $V_n$ with the polinomial ring  $k[\alpha_0,\alpha_1,\ldots,  \alpha_n ]. $
\noindent
The group  $SL_2$ acts on  $V_n$  by the rule
$$
(g\,u)(x,y)=u(d\,x-b\,y, -c\,x+a\,y), \phantom{**} g= \Bigl( \begin{array}{ll} 
 a\, b \\
 c\,d
\end{array} \Bigr) \in SL_2.
$$
\noindent
The generating elements  $\Bigl( \begin{array}{ll}  0\, 1 \\ 0\,0 \end{array} \Bigr),$ $\Bigl( \begin{array}{ll}  0\, 1 \\ 0\,0 \end{array} \Bigr)$ of the tangent Lie algebra  $\mathfrak{sl_{2}}$ act on $V_d$ by  derivations  $-y \frac{\partial}{\partial x}$,  $-x \frac{\partial}{\partial y},$ see \cite{VinP},   and   on  $R_n$  by  derivations
$$
\begin{array}{l}
d_1:=\displaystyle \alpha_0\frac{\partial}{\partial \alpha_1}+2\, \alpha_1\frac{\partial}{\partial \alpha_2}+\cdots +n\,\alpha_{n-1}\frac{\partial}{\partial \alpha_n}, \\
 d_2:=\displaystyle n\,\alpha_1\frac{\partial}{\partial \alpha_0}+(n-1)\,\alpha_2\frac{\partial}{\partial \alpha_1}+\cdots +\alpha_{n}\frac{\partial}{\partial \alpha_{n-1}}.
\end{array}
$$
It follows that the invariant ring  $R_n^{SL_2}$ coincides with a ring of polynomial solutions of the following first order PDE system, see \cite{Hilb}, \cite{Gle}:
$$
\left\{
\begin{array}{r}
\displaystyle \alpha_0\frac{\partial u}{\partial \alpha_1}+2\, \alpha_1\frac{\partial u}{\partial \alpha_2}+\cdots +n\,\alpha_{n-1}\frac{\partial u}{\partial \alpha_n}=0, \\
\displaystyle n\,\alpha_1\frac{\partial u}{\partial \alpha_0}+(n-1)\,\alpha_2\frac{\partial u}{\partial \alpha_1}+\cdots +\alpha_{n}\frac{\partial u}{\partial \alpha_{n-1}}=0,
\end{array}
\right.,
\eqno (*)
$$
i.e. $R_n^{SL_2}=k[\alpha_0,\alpha_1,\ldots,  \alpha_n ] ^{d_1}  \bigcap k[\alpha_0,\alpha_1,\ldots,  \alpha_n ] ^{d_2} ,$ where $u \in k[\alpha_0,\alpha_1,\ldots, \alpha_n], $  and  $$k[\alpha_0,\alpha_1,\ldots,  \alpha_n ] ^{d_i}:=\{ f \in k[\alpha_0,\alpha_1,\ldots,  \alpha_n ]| d_i(f)=0 \}, i=1,2.
$$  
The ring of invariants and theirs syzygies  was a major object of research  in  classical invariant theory of the 19th century.
For  $n\leq 6$ the invariant ring was described by Gordan \cite{Gor}. The case  $n{=}8$  was considered by Shioda   \cite{Shio}.
\noindent
Allmost all known invariants were found in implicit way by the so-called symbolic method. In this case every invariant is represented as an action of an   invariant differential operator   applied to  covariants ($\Omega$ process).
For $n \leq 5$ Fa$\acute {\rm a}$ de Bruno \cite{Bruno}  and Sylvester \cite{Sylv} have calculated an explicit way  a minimal generating set of the  $R_n^{SL_2}.$
 The explicit  form  is highly unwieldy. For example, one invariant of   the generating set for the ring  $R_5^{SL_2}$ is a polynomial of degree 18, consisting of 848 terms. 

 For $n=5,$  Sylvester also   found the explicit form of  single syzygy of degree 36  between the four  invariant  of the generating set. For $n=6,$ there exists a unigue  syzygy and in the case $n=8$ the 9 fundamental invariants are related by 5 syzygies, see \cite{Pop}, \cite{Shio}.

\noindent
It is the aim of this paper to  reduce the system $(*)$ to one equation and try to use it for calculation of invariants for small $n$.
For each  equations of the system  $(*)$ one may easily  find a fundamental system of  solutions in the qoutiend  field $k(\alpha_0,\alpha_1,\ldots, \alpha_n) .$  Having the fundamental system for one derivation  we may discover an action of another  derivation on the fundamental system .  In this way we succeed to reduce the system to single equation. By using the equation together with an additional information about invariants such as theirs degree  and  number of  generators one may calculate an explicit form for invariants and their syzygies. In this paper we do so for  $n \leq 6$ and $n=8.$  

\section{Differential eqation of  invariants}

From now, let us change the variable set   $\alpha_0,\alpha_1,\ldots,  \alpha_n$ on    $x_0,x_1,\ldots, x_n.$  Denote by   $k[X]$ the ring  $k[x_0,x_1,\ldots, x_n],$ and by  $k(X)$ denote a quotient field of the ring  $k[X].$ The derivations  $d_1,d_2$ one may extend from   $k[X]$ to  $k(X)$  preserving the same notations $d_1,d_2.$  It is evident that 
 $k[X]^{d_i}=k(X)^{d_i} \bigcap k[X].$ The derivations   $d_1$ is  localy nilpotent on  $k[X],$ moreover $d_1(\lambda)=-1,$  where $\lambda=-\frac{x_1}{x_0}.$ Therefore  (see for example \cite{Ess}, Proposition 1.3.21)  for the derivation  $d_1$ we can get a description of the ring   $k(X)^{d_1},$  namely 
$$
k(X)^{d_1}=k(\sigma(x_0),\sigma(x_1),\dots ,\sigma(x_n)),
$$
\noindent
where  $\sigma: k[X] \to k[X]^{d_1} $ is a ring homomorfism defined by 
$$
\sigma(a)=\sum_{i=0}^{\infty} d_1^{\,i}(a) \frac{\lambda^i}{i!}.
$$
\noindent
It is well known (\cite{Wri}, Proposition 2.1)  that in this case $k(X)=k(X)^{d_1}[\lambda]$ and $\lambda$ are algebraically independent over $k(X)^{d_1},$ therefore $k(X)$ is a polynomial ring in one variable over $k(\sigma(x_0),\sigma(x_2),\dots ,\sigma(x_n)).$ This  fact allows us to define an action of the derivation $d_2$ on $k(X)$ in new coordinates $\sigma(x_0),\sigma(x_2),\dots ,\sigma(x_n),\lambda.$  Denote by $d$ the derivation
$$
d:=d_2(\lambda)\,\frac{\partial}{\partial \lambda}+d_2(\sigma(x_0))\,\frac{\partial}{\partial\sigma(x_0)}+\ldots +d_2(\sigma(x_n))\,\frac{\partial}{\partial \sigma(x_n)}.
$$ 
Since $k(X)=k(\sigma(x_0),\sigma(x_2),\dots ,\sigma(x_n))[\lambda]$ we are always able to express $d_2(\lambda),$ $d_2(\sigma(x_i))$ in terms of $\lambda, \sigma(x_0),\sigma(x_2),\dots ,\sigma(x_n),$ so the derivation $d$ is well defined.

It is clear that the ring $k[X]^d$ coincides with a ring  polynomial solution of the following differential equation
$$
d_2(\lambda)\,\frac{\partial u}{\partial \lambda}+d_2(\sigma(x_0))\,\frac{\partial u}{\partial \sigma(x_0)}+\ldots +d_2(\sigma(x_n))\,\frac{\partial u}{\partial \sigma(x_n)}=0, u\in k[\sigma(x_0),\sigma(x_1),\dots ,\sigma(x_n)],
\eqno (**)
$$
\noindent
and we get  $R_n^{SL_2}=k[X]^d \bigcap k[X].$

\bigskip
\noindent
For  example, let us consider the case  $n=2$. The derivations  $d_1$ and  $d_2$ has form 
$$
\begin{array}{l}
d_1=\displaystyle x_0\,\frac{\partial}{\partial x_1}+2\,x_1 \frac{\partial}{\partial x_2},\\
d_2=2\,\displaystyle x_1\,\frac{\partial}{\partial x_0}+x_2 \frac{\partial}{\partial x_1}
\end{array}
$$
\noindent
Since 
$$
\begin{array}{l}
\sigma(x_0)=\displaystyle x_0,\\
\sigma(x_1)=\displaystyle x_1+x_0\, \lambda=x_1+x_0\,(-\frac{x_1}{x_0})=0,\\
u_2:=\sigma(x_2)=\displaystyle x_2+2\,x_1\,\lambda +2\,x_0\,\frac{\lambda^2}{2!}=x_2-x_0\,\lambda^2,
\end{array}
$$
\noindent
then  $k(X)^{d_1}=k(x_0,u_2).$
\noindent
Taking into account  $x_2=u_2+x_0\,\lambda^2$ we obtain 
$$
\begin{array}{l}
d_2(x_0)=\displaystyle 2\,x_1=-2 x_0\,\lambda\\
d_2(\lambda)=\displaystyle d_2(-\frac{x_1}{x_0})=\lambda^2-\frac{u_2}{x_0},\\
d_2(u_2):=\displaystyle d_2\Bigl(x_2-x_0\,\lambda^2\Bigr)=2\,\lambda u_2,
\end{array}
$$
Therefore, the derivation  $d$ has form  
$$
d= -2 x_0\,\lambda\frac{\partial}{\partial x_0}+\left(\lambda^2-\frac{u_2}{x_0}\right) \frac{\partial}{\partial \lambda}+2\,x_0 u_2 \frac{\partial}{\partial u_2},
$$
and the corresponding  differential equation  $(**)$  after  multiplying by $x_0$ turns into

$$
-2 x_0^2\,\lambda \frac{\partial u}{\partial x_0}+(x_0 \lambda^2-u_2) \frac{\partial u}{\partial \lambda}+2\,x_0^2 u_2 \frac{\partial u}{\partial u_2}=0.
$$ 
It is well known   that  a first order PDE in three variables has a fundamental system which consists of  no more than 2 solutions.  For our equation it can easily be checked that  the fundamental system, it is  $x_0\,u_2$ and  $u_2+x_0 \lambda^2,$ hence  $k(x_0,u_2)[\lambda]^d=k(x_0\,u_2,u_2+x_0 \lambda^2)$. Further, the ring $k(x_0\,u_2)$ obviously is the intersection of rings   $k(x_0\,u_2,u_2+x_0 \lambda^2)$ and $k(x_0,u_2).$ Since $x_0\,u_2=x_0\,x_2-x_1^2$ already belongs to $ k[X]$ it follows that  the invariant ring  $R_2^{SL_2}$  is generated by single invariant $x_0\,x_2-x_1^2 $.

Consider the case of arbitrary  $n$.  We have  
$$
k(X)^{d_1}=k(x_0,u_2,u_3,\ldots,u_n),
$$
where
$$
u_i:=\sigma(x_i)=\sum_{k=0}^i {i \choose k}\,x_{i-k}\,\lambda^k.
$$
First of all we will express variables   $\lambda, x_2,\ldots,x_n$ through   $u_2,\ldots, u_n$. Denote by $B_i$ the sum $i x_1\lambda^{i-1}+x_0 \lambda^i=-(i-1) x_0 \lambda^{i}$ of the last two terms of  $u_i$. In particulary we obtain 
$$
\begin{array}{l}
u_2=\displaystyle x_2+B_2,\\
u_3=x_3+3\,x_2 \lambda +B_3,\\
u_4=\displaystyle x_4+4\,x_3\lambda+6\,x_2\lambda^2+B_4,
\end{array}
$$
Since
$$
\begin{array}{l}
x_2=\displaystyle u_2-B_2,\\
x_3=u_3-3\,u_2\lambda -(B_3-3\,B_2 \lambda),\\
x_4=u_4-4 u_3 \lambda+6 u_2 \lambda^2-(B_4-4 B_3 \lambda+6 B_2 \lambda^2),
\end{array}
$$
and for arbitrary  $i:$ 
$$
x_i=\sum_{k=0}^{i-2}(-1)^k {i\choose k}u_{i-k}\lambda^k-\Bigl(\sum_{k=0}^{i-2}(-1)^k {i\choose k}B_{i-k}\lambda^k  \Bigr).
$$
Taking into account 
$$
\sum_{k=0}^{i-2}(-1)^k {i\choose k}B_{i-k}\lambda^k=\sum_{k=0}^{i-2} (-1)^{k+1} (i{-}k{-}1) {i\choose k} x_0 \lambda^{i}=
$$
$$
 =x_0 \lambda^{i} \sum_{k=0}^{i-2} (-1)^{k+1} (i-k-1) {i\choose k}=(-1)^{i+1} x_0 \lambda^{i},
$$
we get required form for  $x_i:$
$$
x_i=\sum_{k=0}^{i-2}(-1)^k {i\choose k}u_{i-k}\lambda^k+(-1)^{i} x_0 \lambda^{i}.
$$

To wtite the equation  $(**)$ we need to find an explicit expression of   $d_2(u_i)$ in terms of $\lambda, x_2,\ldots,x_n$ through   $u_2,\ldots, u_n$.
 By direct calculation we obtain

$$
\begin{array}{l}
d_2(\lambda)=\displaystyle \lambda^2-(n-1) \frac{u_2}{x_0},\\
d_2(u_2)=\displaystyle(n-2)u_3-(n-4)u_2 \lambda,\\
d_2(u_3)=\displaystyle (n-3)u_4 -(n-6)u_3 \lambda-3(n-1)\frac{u_2^2}{x_0}.
\end{array}
$$
In the general case,  for  $u_i, i>3$ we have: 
$$
d_2(u_i)=\sum_{k=0}^{i-2} {i \choose k}\,d_2(x_{i-k})\,\lambda^k+\left(\sum_{k=1}^{i-2} k x_{i-k}  {i \choose k} \lambda^{k-1} \right) d_2(\lambda) +d_2(B_i).
$$
Let us calculate each sum separately
$$
\begin{array}{l}
\mbox{    }\displaystyle \sum_{k=0}^{i-2} {i \choose k}\,d_2(x_{i-k})\,\lambda^k=\sum_{k=0}^{i-2} {i \choose k}\,(n{-}(i{-}k))\,x_{i-k+1}\lambda^k=\\

\displaystyle=\sum_{k=0}^{i-2} {i \choose k}\,(n{-}(i{-}k))\,\lambda^k \left( \sum_{s=0}^{i-k-1}(-1)^s {i-k+1\choose s} u_{i-k-s+1}\lambda^s +(-1)^{i-k+1} x_0 \lambda^{i-k+1}\right)=\\

\displaystyle=\sum_{k=0}^{i-2}\sum_{s=0}^{i-k-1} (-1)^{s}(n{-}(i{-}k)) {i \choose k}{i-k+1\choose s} u_{i-k-s+1}\lambda^{s+k}+x_0 \lambda^{i+1} T_i=\\

\displaystyle=\sum_{p=2}^{i+1} u_{p}\lambda^{i+1-p} \sum_{s+k=i+1-p} (-1)^s (n{-}(i{-}k)) {i \choose k}{i-k+1\choose s} +x_0 \lambda^{i+1}T_i=\\

\displaystyle=u_2 \lambda^{i-1}S_2+\sum_{p=3}^{i+1}u_p \lambda^{i-p+1}\,S_p+x_0 \lambda^{i+1} T_i=\\
\end{array}
$$
here
$$
\begin{array}{l}
T_i:=\displaystyle \sum_{k=0}^{i-2} (-1)^{i-k+1} (n{-}(i{-}k)){i \choose k},\\
S_2:=\displaystyle \sum_{k=3}^{i+1} (-1)^{k-2}(n-(k-1)) {k \choose 2}{i\choose k-1},\\
S_p:=\displaystyle \sum_{k=p}^{i+1} (-1)^{k-p}(n-(k-1)) {k \choose p}{i\choose k-1}, p>2 \\
\end{array}
$$
\begin{lm} The following equalities  are hold:
$$
T_i=n+i-n\, i, i>1.
$$
and for  $i>3$
$$
S_2=-(n-1)i
$$
$$
S_p=\left\{ \begin{array}{l}
n-i, \mbox{\rm {for} } p=i+1,\\
2\,i-n, \mbox{\rm {for} } p=i,\\
-i,  \mbox{\rm {for} } p=i-1,\\
0, \mbox{\rm {for} }2<p<{i-1}, \\
-(n-1)i,\mbox{\rm {for} }p=2.
\end{array}
 \right.
$$
\end{lm}
\begin{proof}
 Using the  binomial identity 
$$
\sum_{k=0}^{i}(-1)^{k} {i\choose k}=\sum_{k=1}^{i-1}(-1)^{k-1}{i-1 \choose k-1}=0,
$$
we get
$$
T_i-n+(n-1)i= \sum_{k=0}^{i} (-1)^{i-k+1} (n{-}(i{-}k)){i \choose k}=(n-i)(-1)^{i+1}\sum_{k=0}^{i}(-1)^{-k}{i \choose k} +
$$
$$
+(-1)^{i+1}\sum_{k=0}^{i}(-1)^{-k}k {i \choose k}=0+(-1)^{i+1}\sum_{k=1}^{i-1}(-1)^{k-1}{i-1 \choose k-1}=0.
$$
The equalities for  $S_p$  are derived from the  following ortogonal relation  (\cite{Rio} )  
$$
\delta_{m,n}=\sum_{k=m}^{n} (-1)^k {k\choose m}{n\choose k}
$$
in the same way.
\end{proof}
Hence we reduced the first sum to the form
$$
 (n-i)u_{i+1}-(n-2 i)u_i \lambda-i u_{i-1} \lambda^2-i (n-1)u_2 \lambda^{i-1} +(n+i-n\,i) x_0\lambda^{i+1}.
$$
Let us calculate the second sum in the expression of  $d_2(u_i):$ 
$$
\sum_{k=1}^{i-2} x_{i-k} k {i \choose k} \lambda^{k-1} d_2(\lambda) =\sum_{k=1}^{i-2} x_{i-k} i {i-1 \choose k-1} \lambda^{k-1} d_2(\lambda)=i (u_{i-1}-B_{i-1}) d_2(\lambda).
$$
and   $d_2(B_i):$
$$
\begin{array}{l}
d_2(B_i)=\displaystyle d_2(-(i-1) x_0 \lambda^2)= (i-1)\lambda^{i-1}((n-i)x_0\lambda^2+i(n-1)u_2).
\end{array}
$$
Thus after all  simplifications for $i>3$  we get 
$$
d_2(u_i)=(n-i)u_{i+1}-(n-2i)u_i \lambda-i(n-1) \frac{u_2 u_{i-1}}{x_0}.
$$
Consequently the derivation $d$ acts on   $k(x_0,\lambda,u_2,\dots,u_n)$ by the rule
$$
\begin{array}{l}
d(x_0)=-n x_0 \lambda,\\
d(\lambda)=\displaystyle \lambda^2-(n-1) \frac{u_2}{x_0},\\
d(u_2)=\displaystyle(n-2)u_3-(n-4)u_2 \lambda,\\
d(u_i)=\displaystyle (n-i)u_{i+1}-(n-2i)u_i \lambda-i(n-1) \frac{u_2 u_{i-1}}{x_0}, \mbox{for   }  i>2.
\end{array}
$$

Finally, we obtain

\begin{te}
The invariant ring  $R_n^{SL_2}$ coincides with a ring of polynomial solutions of the following first order PDE
$$
\begin{array}{l}
\displaystyle-n x_0^2 \lambda \frac{\partial u}{\partial x_0}+(x_0 \lambda^2{-}(n{-}1)u_2) \frac{\partial u}{\partial \lambda }+((n-2) u_3 x_0-(n-4)u_2 x_0 \lambda) \frac{\partial u}{\partial u_2 }+\\
\displaystyle+\sum_{i=4}^{n}((n-i)u_{i+1} x_0-(n-2i)u_i x_0 \lambda-i(n-1) u_2 u_{i-1}) \frac{\partial u}{\partial u_i }=0,
\end{array}
$$
where, 
$$
\begin{array}{l}
\lambda=\displaystyle \frac{x_1}{x_0},\\
u_i=\displaystyle \sum_{k=0}^i {i \choose k}\,x_{i-k}\,\lambda^k,
\end{array}
$$
and  $u\in k[X] \bigcap k[x_0,u_2,\ldots,u_n].$ 
\end{te}

\section{Solving of the invariant equation}
Let us introduce on $k[x_0,u_2,\dots, u_n]$ three additional derivations $\hat d_1$,$\hat d_2$ and $e$ as follows
$$
\begin{array}{l}
\hat d_1(u_i)=i\,u_{i-1}, \hat d_1(u_2)=0, \hat d_1(x_0)=0,\\
\hat d_2(u_i)=(n-i)\,u_{i+1}, \hat d_2(x_0)=0,\\
e(u_i):=[\hat d_1,\hat d_2](u_i)=(n-2\,i) u_i, e(x_0)=n\,x_0.
\end{array}
$$
Then,  in terms of these derivations, one may rewrite the derivation $d$ in the form
$$
d=x_0\,\hat d_2-x_0 \lambda e-(n-1) u_2 \hat d_1. \eqno (***)
$$
 For the derivation  $\,e,$ every monomial $x_0^{\alpha_1}\, u_2^{\alpha_2}\,\cdots u_n^{\alpha_n}$  is an eigenvector with the eigenvalue \begin{equation} \omega( x_0^{\alpha_0}\, u_2^{\alpha_2}\,\cdots u_n^{\alpha_n}) =n\,(\sum_i\,\alpha_i)-2\,(\alpha_2+2\,\alpha_2+\cdots +n\,\alpha_n).  \notag\end{equation} 
 A homogeneous polynomial is called isobaric if  the sum $(2\,\alpha_2+3\,\alpha_3+\cdots +n\,\alpha_n)$ has an equal value on all  monomials of the polynomial.  The value  is called   $u$-weight of the polynomial  $z$ and it is denoted by $\omega_u(z)$. Hence, if $z$ is isobaric polynomial then  we have $e(z)=(n \deg(z)-2\,\omega_u(z))\,z.$
\noindent
  Put
 $$I_n:=\{ z \in k[x_0,u_2,\ldots ,u_n],\mbox{   } n\,\deg(z)-2\,\omega_u(z)=0 \}.$$
 It is clear that on  isobaric polynomials the functions 
  $\deg, \omega_u:k[x_0,u_2,\ldots ,u_n] \to \mathbb{Z} $ are  additive. Thus $I_n$ form a  subring of the ring  $k[x_0,u_2,\ldots ,u_n].$ 
\noindent
The following theorem describes  solutions of the equation for the invariants.
\begin{te}
$R_n^{SL_2}=I_n \cap k[X].$
\end{te}
\begin{proof}
Suppose that a polynomial  $u \in R_n^{SL_2}.$ Then  $u \in  k[x_0,u_2,\ldots ,u_n]$ and it is follow  $\hat d_i(u) \in  k[x_0,u_2,\ldots ,u_n],i=1,2.$  
Therefore the equality  
$$
d(u)=x_0\,\hat d_2(u)-x_0 \lambda e(u)-(n-1) u_2 \hat d_1(u)=0,
$$
is possible only in the case if a coefficient by  $\lambda$ is equal to zero, i.e. $e(u)=0.$ Therefore  $n\,\deg(u) -\omega_u(u)=0$ and $u \in I_n \cap k[X]$. So we get   $R_n^{SL_2} \subset I_n \cap k[X].$

Now suppose  $u \in I_n \cap k[X] $.  Let us show that  $d_2(u)=0.$ 
Define a  $k$-linear, multiplicative map  $\varphi:k[X] \to k[x_0,u_2,\ldots ,u_n]$ by the rule $\varphi(x_0)=x_0,$ $\varphi(x_1)=0$ and  $\varphi(x_i)=u_i$  for  $ 2 \leq i \leq n.$
\begin{lm}
 If  $z \in k[x_0,u_2,\ldots ,u_n],$ then  $\varphi(z)=z$.
\end{lm}
\begin{proof}
Using the expansion  $$x_i=\sum_{k=0}^{i-2}(-1)^k {i\choose k}u_{i-k}\lambda^k+(-1)^{i} t \lambda^{i},$$  we may every  $x_i$ expres in  form  
 $x_i=u_i+\lambda\,u'_i=\varphi(x_i)+\lambda\,u'_i$  for some $u_i'.$ Multiplicativity of  maps  $\varphi$  and  $x_1=-\lambda x_0=\varphi(x_1) -\lambda x_0,$  implies the existence of the representation  $z=\varphi(z)+\lambda z'$ for arbitrary polynomial  $z$ and for some polynomial  $z'.$ Since  $z \in k[x_0,u_2,\ldots ,u_n]$  and  $\lambda z' \notin k[x_0,u_2,\ldots ,u_n]$ it is follows  $z'=0,$ and therefore   $\varphi(z)=z.$

Consider, for example, the polynomial  $z:=x_4 x_0-4\,x_1\,x_3+3\,x_2^2$. Substituting 
$$
\begin{array}{l}
x_1 = -\lambda x_0,\\
x_2 = u_2+x_0 \lambda^2,\\
x_3 = u_3-3\,u_2\,\lambda-x_0\, \lambda^3,\\
x_4 =u_4+ 6 u_2 \,\lambda^2-4\,u_3\,\lambda+x_4\,\lambda^4,
\end{array}
$$
in $z$, after simplification we obtain 
$z=x_0\,u_4+3\,u_2^2$. On the other hand $$\varphi(z)=\varphi(x_4\,x_0-4\,x_1\,x_3+3\,x_2^2)=x_0\,u_4+3\,u_2^2.$$ Thus  $\varphi(z)=z.$
\end{proof}

Recall now that a weight  $\omega(z)$ of homogeneous isobaric invariant  $z \in k[X]^{SL_2}$ is called the value $n\,\deg(z)-\omega(z),$ where $\omega$ is a integral function which takes each monomial   $x_0^{\alpha_0}\, x_1^{\alpha_1}\,\cdots x_n^{\alpha_n}$ to  $\alpha_1+2\,\alpha_2+\cdots +n\,\alpha_n.$  From previous lemma it is follows that $\omega(z)=\omega_u(z)$ if only    $z \in I_n \cap k[X].$ It is well known  (Hilbert, \cite{Hilb} , p.38) that for an isobaric polynomial  $z,$ the conditions   $d_1(z)=0$ and  $n\deg(z)=\omega(z)$ follow  $d_2(z)=0.$ Thus each  polynomial of  $ I_n \cap k[X]$ is invariant of  $SL_2.$
\end{proof}

\section{Algorithm }

By using the results of theorems 1, 2, in the cases  $n \leq 6$  and  $n=8$, one may develope an  effective algorithm for calculations of minimal generating sets of the invariant ring $R_n^{SL_2}$ and theirs syzygies. The main algorithm consist of several subsidiary algorithms. Note that we use the  information about number of invariants, number of syzygies and theirs degree.

Let us go to the description of the algorithms.
\bigskip
\begin{center}  {\bf Algorithm (main)}  ${\rm MINGENSET(n,r,D)}.$\end{center}
\noindent
{\bf Input:} $n\,$ is degree of an binary form; $r$ is a number of homogeneous invariants which forms a minimal generating set of the invariant ring  $R_n^{SL_2}$ and  $D:=\{s_1,s_2,\ldots, s_r \}, s_i\leq s_{i+1}$ is the set of their degries  \\
{\bf Output:} Minimal generating set of the invariant ring $R_n^{SL_2}.$\\
\noindent
{\bf begin}\\
\noindent
$S:=\{\emptyset\};$\\
\noindent
{\bf for i from $s_1$  to $s_r$ do } \\
\noindent
$I={\rm INVARIANTS(n,i)};$\\
\noindent
{\bf for k from $1$  to $nops(I)$ do } \\
{\bf if ${\rm MEMBER}(S,I_k)$ then $S:=S$ union $ \{I_k\}$ end if;}\\
\noindent
{\bf end do;\\
\noindent
end do;\\
\noindent
return S;\\
\noindent
end. }\\
\bigskip
\begin{center}  {\bf Algorithm }  ${\rm INVARIANTS}(n,d).$ \end{center}
\noindent
{\bf Input:} $n $ is the degree of  binary form;$d\,$ is a degree of an invariant.\\
{\bf Output:} The set of linearly independed homogeneous invariants of degree  $n.$\\
\noindent
{\bf begin}\\
\noindent
$P:={\rm POWERS}(n,d);$\\
\noindent
$F:={\rm GPOL}(P);$\\
\noindent
$F1:={\rm DERPOL}(F);$ \\
\noindent
$F2:={\rm COEFFS}(F1);$ \\
\noindent
$F3:={\rm SYSTEM}(F2);$\\
\noindent
$I:={\rm SOLVSUBS}(F3,F);$ \\
\noindent
{\bf return $I;$ \\
\noindent
end.}

\bigskip
\begin{center}  {\bf Algorithm}  ${\rm MEMBER}(S,F).$ \end{center}
\noindent
{\bf Input:} $S\,$ is a set of  homogeneous polynomials $\{f_1,f_2,\ldots , f_m \};$ $f\,$ a homogeneous polynomial.\\
{\bf Output:} ${\rm TRUE}$  if  $F \notin k[f_1,f_2,\ldots , f_m ]$ and   ${\rm FALSE},$ if $F \in k[s_1,s_2,\ldots , s_m ].$ \\
\noindent
{\bf begin}\\
\noindent
$P:={\rm GRAD}(S,f);$\\
\noindent
$F:={\rm GPOL2}(P);$\\
\noindent
$F1:={\rm SUBS}(S,F)$\\
\noindent
$F2:={\rm COEFFS}(F1);$ \\
\noindent
$F3:={\rm SYSTEM}(F2);$\\
\noindent
$S:={\rm ISSOLVABLE}(F3);$\\
\noindent
{\bf return $S;$\\
end.}

\bigskip
\begin{center}  {\bf Algorithm}  ${\rm SYZYGIES}(S,r,D).$ \end{center}
\noindent
{\bf Input:} $S\,$ is a set of homogeneous polynomials  $\{f_1,f_2,\ldots , f_m \};$ $r\,$ is a number of syzygies and   $D:=\{s_1,s_2,\ldots, s_r \}$ is a set of their degries.\\
{\bf Output:}The set of $r$ syzygies  of the set  $S$ .  \\
\noindent
{\bf begin}\\
\noindent
$S:=\{\emptyset\};$\\
\noindent
{\bf for i from $1$  to $r$ do } \\
\noindent
$F:={\rm POWERS2}(S,i);$\\
\noindent
$F1:={\rm GPOL2}(P)$\\
\noindent
$F2:={\rm COEFFS}(F1)$\\
\noindent
$F3:={\rm SYSTEM}(F)$\\
\noindent
$F4:={\rm SOLVSYSTEM}(R,F)$\\
\noindent
$S:=S$    {\bf union} $ \{F4 \};$\\
\noindent
{\bf end do;\\
\noindent
return S;\\
\noindent
end. }\\

\noindent
\bigskip
\begin{center}  {\bf Algorithm}  ${\rm POWERS}(n,d).$ \end{center}
{\bf Input:} $n\,$ is binary form degree; $d\,$ is a degree of an invariant.\\
{\bf Output:} A set of solutions of the equation system 
$$
\left \{
\begin{array}{l}
\alpha_0+\alpha_2+\ldots +\alpha_n=d,\\
\displaystyle 2 \alpha_2+3 \alpha_3 +\ldots +n \alpha_n= \frac{n\,d}{2},
\end{array}
\right.
\alpha_i \in \mathbb Z_{+}.
$$
\noindent
\bigskip
\begin{center}  {\bf Algorithm }  ${\rm POWERS2}(S,d).$ \end{center}
{\bf Input:} $S\,$ is set of homogeneous polynomials $\{f_1,f_2,\ldots , f_m \};$ $d\,$ is a degree of their sysygy.\\
{\bf Output:} $P\,$ is a set of solutions of the equation 
$$
\alpha_1 \deg(f_1)+\alpha_2 \deg(f_2)+\ldots +\alpha_m \deg(f_m)=d , \alpha_i \in \mathbb Z_{+}.
$$

\bigskip
\begin{center}  {\bf Algorithm}  ${\rm GRAD}(S,f).$ \end{center}
{\bf Input:} $S\,$ is a set of homogeneous polynomials $\{f_1,f_2,\ldots , f_m \};$ $f\,$ is a homogeneous polynomial.\\
{\bf Output:} $P\,$ is set of solutions of the equation   
$$
\left \{
\begin{array}{l}
\alpha_1 \deg(f_1)+\alpha_2 \deg(f_2)+\ldots +\alpha_m \deg(f_m)=\deg(f) ,\\
\alpha_1 \omega(f_1)+  \alpha_2 \omega(f_2) +\ldots + \alpha_m \omega(f_m)=\omega(f). 
\end{array}
\right.
\alpha_i \in \mathbb Z_{+}.
$$

\bigskip
\begin{center}  {\bf Algorithm}  ${\rm GPOL}(P).$ \end{center}
{\bf Input:} $P\,$ is a set of $n$-tuples $\{(\alpha_0,\alpha_2,\ldots, \alpha_n )\}.$\\
{\bf Output:} $F\,$ is a "general" polinomial  $$ F:=\sum_{\alpha \in P} \beta_{\alpha} x_0^{\alpha_0} u_2^{\alpha_2}\cdots u_n^{\alpha_n}.$$

\bigskip
\begin{center}  {\bf Algorithm}  ${\rm GPOL2}(P).$ \end{center}
{\bf Input:} $P\,$  is a set of $m$-tuples $\{(\alpha_1,\alpha_2,\ldots, \alpha_m )\}.$\\
{\bf Output:} $F\,$  is a "general" polinomial  $$ F:=\sum_{\alpha \in P} \beta_{\alpha} f_1^{\alpha_1} f_2^{\alpha_2}\cdots f_m^{\alpha_m}.$$

\bigskip
\begin{center}  {\bf Algorithm}  ${\rm DERPOL}(F).$ \end{center}
{\bf Input:} $F\,$ is a polynomial of  $F \in k[x_0,u_2,\ldots, u_n].$ \\
{\bf Output:} A result of the action of the derivation $d:=x_0\,\hat d_2-(n-1) u_2 \hat d_1$ applied to the polynomial  $F.$

\bigskip
\begin{center}  {\bf Algorithm}  ${\rm SUBS}(S,F).$ \end{center}
{\bf Input:}$S\,$   is a set of homogeneous polynomials  $\{f_1,f_2,\ldots , f_m \};$ $F\,$ is a polynomial  of  \\ $F \in k[f_1,f_2,\ldots , f_m].$ \\
{\bf Output:} A result of substituting of $S$ in the polynomial  $F$.  

\bigskip
\begin{center}  {\bf Algorithm}  ${\rm COEFFS}(F).$ \end{center}
{\bf Input:} $F\,$ is a polynomial of  $k[x_0,u_2,\ldots, u_n].$ \\
{\bf Output:} A set of coefficients of the polynomial $F.$

\bigskip
\begin{center}  {\bf Algorithm}  ${\rm SYSTEM}(F).$ \end{center}
{\bf Input:} $F\,$ a polynomial  $F \in  k[x_0,u_2,\ldots, u_n].$ \\
{\bf Output:} A system of linear homogeneous equations for the indeterminates $ \beta_1,\beta_2 ,\ldots $. We get the system  by using the condition $F \equiv 0.$

\bigskip
\begin{center}  {\bf Algorithm}  ${\rm SOLVSYSTEM}(R,F).$ \end{center}
{\bf Input:} $R\,$ is a system of linear homogeneous equations; $F\,$  is a  "general"  polynomial of  $k[f_1,f_2,\ldots, f_m].$ \\
{\bf Output:} A result of substituting of solutions of the sysytem  $R$ in the polynomial  $F.$

\bigskip
\begin{center}  {\bf Algorithm}  ${\rm SOLVSUBS}(R,F).$ \end{center}
{\bf Input:} $R\,$ is a system of linear homogeneous equations for the indeterminates $ \beta $;\\ $F\in  k[x_0,u_2,\ldots, u_n].$ \\
{\bf Output:} A  vector space basis of  solutions of the equation $d(F)=0.$

\bigskip
\begin{center}  {\bf Algorithm}  ${\rm ISSOLVABLE}(R).$ \end{center}
{\bf Input:} $R\,$ is a system of linear homogeneous equations for the indeterminates $ \beta_1,\beta_2 ,\ldots $ \\
{\bf Output:} ${\rm TRUE}$ if the system $R$ has no solutions  and  ${\rm FALSE}$
if  $R$   has solutions. \\

\section{Examples}
Below one may find   differential equations  and syzygies   of the invariant ring  $R_n^{SL_2}$ for  $n\leq 6 $ and $n=8.$ Also here is plased a minimal generating system of the invariants rings in the case $n<5.$ 
All calculations were done with Maple according to the above algorithms.
 
\bigskip
\noindent
${\bf n=3.}$ 

The differential equation for the invariant $u$,  $u\in k[x_0,u_2,u_3]$  is as follow:
$$
 - 3\,x_0^{2}\,\lambda \,{\frac {\partial }{\partial x_0}}\,u
 +  ({u_{3}}\,x_0{+}{u_{2}}\,x_0\,\lambda )\,({\frac {
\partial }{\partial {u_{2}}}}\,u){+}(3\,{u_{3}}\,x_0{-}6\,{u_{2}}^{
2})\,({\frac {\partial }{\partial {u_{3}}}}\,u)=0,
$$
or 
$$
d(u)=x_0\,\hat d_2(u)-2 u_2 \hat d_1(u)=0.
$$
The invariant ring  $R_3^{SL_3}$ generated by one invariant $f_4$ of degree four. The subscript in $f_i$ means a degree of the polynomial $f_i.$The weight  of $f_4$  is equal  $\displaystyle \frac{3\cdot 4}{2}=6.$ The system of equations
$$
\left \{
\begin{array}{l}
\alpha_0+\alpha_2 +\alpha_3=4,\\
\displaystyle 2 \alpha_2+3 \alpha_3 =6, 
\end{array}
\right.
$$
in $\mathbb Z_{+}\,$ has only the following two solutions $-$ $(1,3,0)$  i $( 2,0,2).$ 
Then we find an invariant in the form  $f_4=\beta_1 x_0 u_2^3+\beta_2  x_0^2 u_3^2. $ We have
$$
d(f_4)=x_0\,\hat d_2(f_4)-2 u_2 \hat d_1(f_4) =x_0(\beta_1 3 x_0 u_2^2 u_3)-2 u_2 (\beta_2 6 x_0^2 u_2 u_3)=(3 \beta_1 -12 \beta_2) x_0^2 u_2^2 u_3=0.
$$
This implyies $3 \beta_1 -12\beta_2=0$, or  $\beta_1=4 \beta_2.$  Thus $f_4=4 x_0 u_2^3+x_0^2 u_3^2$  and 
$$
R_3^{SL_2}=k[f_4].
$$
Going back to indeterminates  $x_0,x_1,x_2,x_3$ we get 
$$
f_4=4\, x_0 x_2^3-3 x_1^2 x_2^2+t^2\,x_3^2-6\,x_0 x_1 x_2 x_3+4\,x_1^3\,x_3.
$$
\noindent
${\bf n=4.}$ 

The differential equation for the invariant $u$,  $u\in k[x_0,u_2,u_3,u_4]$  is as follow:

$$
 - 4\,x_0^{2}\,\lambda \,({\frac {\partial }{\partial t}}\,u)
+ 2\,{u_{3}}\,x_0\,({\frac {\partial }{\partial {u_{2
}}}}\,u) + (x_0\,{u_{4}} + 2\,\lambda \,x_0\,{u_{3}} - 9\,{u_{2}}^{2})
\,({\frac {\partial }{\partial {u_{3}}}}\,u) \\
\mbox{} + 
$$
$$
+(4\,{u_{4}}\,x_0\,\lambda  - 12\,{u_{2}}\,{u_{3}})\,(
{\frac {\partial }{\partial {u_{4}}}}\,u) =0
$$
The invariant ring  $R_4^{SL_2}$ generated by two invariants $f_2, f_3.$ Those invariants one may find in the same way as it were done for the case $n=3.$.

$$
R_4^{SL_2}=k[f_2,f_3], {f_{2}} := t\,{u_{4}} + 3\,{u_{2}}^{2},{f_{3}} := {u_{2}}^{3} - t\,{u_{2}}\,{u_{4}} + t\,{u_{3}}^{2}
$$
\noindent
${\bf n=5.}$ 

The differential equation is as follow 
$$
- 5\,x_0^{2}\,\lambda \,({\frac {\partial }{\partial x_0}}\,u)
 + (3\,{u_{3}}\,x_0 - {u_{2}}\,x_0\,\lambda )\,({\frac {
\partial }{\partial {u_{2}}}}\,u) \\
\mbox{} + (2\,x_0\,{u_{4}} + \lambda \,x_0\,{u_{3}} - 12\,{u_{2}}^{2})
\,({\frac {\partial }{\partial {u_{3}}}}\,u) +
$$
$$
 +({u_{5}}\,x_0 + 3\,{
u_{4}}\,x_0\,\lambda  - 16\,{u_{2}}\,{u_{3}})\,({\frac {\partial }{
\partial {u_{4}}}}\,u) \\
\mbox{} + (5\,{u_{5}}\,x_0\,\lambda  - 20\,{u_{2}}\,{u_{4}})\,(
{\frac {\partial }{\partial {u_{5}}}}\,u)
 =0
$$
The invariant ring  $R_5^{SL_2}$ generated by four  invariants $f_4,f_8,f_{12},f_{18}.$ 

$$
R_5^{SL_2}=k[f_4,f_8,f_{12},f_{18}],
$$
There exists the single sygyzy:
$$
\begin{array}{l}
 1296\,f_{18}^2= - 48\,{f_{12}}^{3} + {f_{4}}^{5}\,{f_{8}}^{2} - 6\,{f_{4}}^{3}\,
{f_{8}}^{3} + 9\,{f_{4}}\,{f_{8}}^{4} - 2\,{f_{4}}^{4}\,{f_{8}}\,
{f_{12}} - 18\,{f_{4}}^{2}\,{f_{8}}^{2}\,{f_{12}} + 72\,{f_{8}}^{
3}\,{f_{12}} +  \\
\mbox{}  + {f_{4}}^{3}\,{f_{12}}^{2}+ 72\,{f_{4}}\,{f_{8}}\,{f_{12}}^{2}
\end{array}
$$
\noindent
${\bf n=6.}$ 

The diffrential equation is as follows
$$
- 6\,x_0^{2}\,\lambda \,{\frac {\partial }{\partial x_0}}\,u + (4\,
{u_{3}}\,x_0 - 2\,{u_{2}}\,x_0\,\lambda )\,{\frac {\partial }{
\partial {u_{2}}}}\,u +
$$
$$
+(3\,{u_{4}}\,x_0 - 15\,{u_{2}}^{2})\,
{\frac {\partial }{\partial {u_{3}}}}\,u+ (2\,{u_{5}}\,x_0 + 2\,{u_{4}}\,x_0\,\lambda  - 20\,{u_{2}}\,
{u_{3}})\,{\frac {\partial }{\partial {u_{4}}}}\,u\\
$$
$$
 + ({u_{6}}\,
x_0 + 4\,{u_{5}}\,x_0\,\lambda  - 25\,{u_{2}}\,{u_{4}})\,{\frac {
\partial }{\partial {u_{5}}}}\,u
\mbox{} + (6\,{u_{6}}\,x_0\,\lambda  - 30\,{u_{2}}\,{u_{5}})\,
{\frac {\partial }{\partial {u_{6}}}}\,u =0
$$
The invariant ring  $R_6^{SL_2}$ generated by five  invariants $f_2,f_4,f_6,f_{10},f_{15}.$ 

$$
R_6^{SL_2}=k[f_2,f_4,f_6,f_{10},f_{15}],
$$
The exists the single sygyzy 
$$
\begin{array}{l}
640722656250000\,{f_{15}}^{2}=4556250000\,{f_{2}}\,{f_{4}}^{3}\,{
f_{6}}\,{f_{10}} - 6378750000000\,{f_{2}}^{2}\,{f_{4}}\,{f_{6}}^{
2}\,{f_{10}} \\
\mbox{} + 89606250000\,{f_{2}}^{5}\,{f_{4}}\,{f_{6}}\,{f_{10}} + 
51637500000\,{f_{2}}^{3}\,{f_{4}}^{2}\,{f_{6}}\,{f_{10}} + 160692
\,{f_{2}}^{11}\,{f_{4}}^{2} \\
\mbox{} + 864\,{f_{2}}^{3}\,{f_{4}}^{6} + 10242\,{f_{2}}^{5}\,{f
_{4}}^{5} + 54820\,{f_{2}}^{7}\,{f_{4}}^{4} + 133195\,{f_{2}}^{9}
\,{f_{4}}^{3} + 94864\,{f_{2}}^{13}\,{f_{4}} \\
\mbox{} + 27\,{f_{2}}\,{f_{4}}^{7} - 20344800\,{f_{2}}^{12}\,{f_{
6}} - 14175\,{f_{4}}^{6}\,{f_{6}} + 6473250000\,{f_{2}}^{9}\,{f_{
6}}^{2} \\
\mbox{} - 817003125000\,{f_{2}}^{6}\,{f_{6}}^{3} + 379687500\,{f
_{4}}^{3}\,{f_{6}}^{3} + 30691406250000\,{f_{2}}^{3}\,{f_{6}}^{4}
 \\
\mbox{} - 158760000\,{f_{2}}^{10}\,{f_{10}} - 1822500\,{f_{4}}^{5
}\,{f_{10}} - 3189375000000\,{f_{2}}^{5}\,{f_{10}}^{2} \\
\mbox{} - 71191406250000\,{f_{6}}^{5} + 2562890625000000\,{f_{10}
}^{3} + 21952\,{f_{2}}^{15} \\
\mbox{} + 14207062500\,{f_{2}}^{7}\,{f_{4}}\,{f_{6}}^{2} - 540675
\,{f_{2}}^{2}\,{f_{4}}^{5}\,{f_{6}} + 9055125000\,{f_{2}}^{5}\,{f
_{4}}^{2}\,{f_{6}}^{2} \\
\end{array}
$$
${\bf n=8.}$ 

The differential equations
$$
\begin{array}{l}
\displaystyle- 8\,x_0^{2}\,\lambda \,({\frac {\partial }{\partial x_0}}\,u) + (6\,
{u_{3}}\,x_0 - 4\,{u_{2}}\,x_0\,\lambda )\,({\frac {\partial }{
\partial {u_{2}}}}\,u) + (5\,x_0\,{u_{4}} - 2\,\lambda \,x_0\,{u_{3}}
 - 21\,{u_{2}}^{2})\,({\frac {\partial }{\partial {u_{3}}}}\,u)
 \\
\displaystyle\mbox{} + (4\,{u_{5}}\,x_0 - 28\,{u_{2}}\,{u_{3}})\,({\frac {
\partial }{\partial {u_{4}}}}\,u) + (3\,{u_{6}}\,x_0 + 2\,{u_{5}}\,
x_0\,\lambda  - 35\,{u_{2}}\,{u_{4}})\,({\frac {\partial }{\partial 
{u_{5}}}}\,u) \\
\displaystyle\mbox{} + (2\,{u_{7}}\,x_0 + 4\,{u_{6}}\,x_0\,\lambda  - 42\,{u_{2}}\,
{u_{5}})\,({\frac {\partial }{\partial {u_{6}}}}\,u) + ({u_{8}}\,
x_0 + 6\,{u_{7}}\,x_0\,\lambda  - 49\,{u_{2}}\,{u_{6}})\,({\frac {
\partial }{\partial {u_{7}}}}\,u) \\
\displaystyle\mbox{} + (8\,{u_{8}}\,x_0\,\lambda  - 56\,{u_{2}}\,{u_{7}})\,(
{\frac {\partial }{\partial {u_{8}}}}\,u)=0
\end{array}
$$
The invariant ring  $R_8^{SL_2}$ generated by nine invariants $f_2,f_3,f_4,f_{5},f_{6},f_7,f_8,f_9,f_{10}.$ 
\noindent
There exists the five sygyzies
$$
\begin{array}{l}
100822400\,{f_{5}}\,{f_{4}}\,{f_{3}}\,{f_{2}}^{2} + 2318131200\,{
f_{6}}\,{f_{5}}\,{f_{3}}\,{f_{2}} + 10838016\,{f_{6}}^{2}\,{f_{4}
} - 94832640\,{f_{9}}\,{f_{7}} \\
\mbox{} - 67436544\,{f_{10}}\,{f_{6}} - 702464\,{f_{10}}\,{f_{2}}
^{3} + 127573232640\,{f_{5}}^{2}\,{f_{2}}^{3} + 336672\,{f_{3}}^{
2}\,{f_{2}}^{5} \\
\mbox{} - 643328\,{f_{6}}\,{f_{2}}^{5} + 1628\,{f_{4}}\,{f_{2}}^{
6} + 4070\,{f_{4}}^{2}\,{f_{2}}^{4} - 1303680\,{f_{8}}\,{f_{2}}^{
4} - 31911936\,{f_{6}}^{2}\,{f_{2}}^{2} \\
\mbox{} - 8996\,{f_{4}}^{3}\,{f_{2}}^{2} - 8956416\,{f_{3}}^{4}\,
{f_{2}}^{2} - 12117504\,{f_{4}}\,{f_{3}}^{4} - 8768289116160\,{f
_{5}}^{2}\,{f_{3}}^{2} \\
\mbox{} + 48470016\,{f_{10}}\,{f_{3}}^{2} - 2808960\,{f_{8}}\,{f
_{4}}^{2} + 12324702781440\,{f_{6}}\,{f_{5}}^{2} + 411936\,{f_{4}
}^{2}\,{f_{3}}^{2}\,{f_{2}} \\
\mbox{} - 1369804800\,{f_{5}}\,{f_{3}}^{3}\,{f_{2}} + 4741632000
\,{f_{9}}\,{f_{5}}\,{f_{2}} - 253848161280\,{f_{5}}^{2}\,{f_{4}}
\,{f_{2}} \\
\mbox{} + 1404928\,{f_{10}}\,{f_{4}}\,{f_{2}} - 88670400\,{f_{5}}
\,{f_{4}}^{2}\,{f_{3}} + 8851046400\,{f_{8}}\,{f_{5}}\,{f_{3}} - 
8467200\,{f_{9}}\,{f_{4}}\,{f_{3}} \\
\mbox{} - 58329600\,{f_{7}}\,{f_{5}}\,{f_{4}} + 71688960\,{f_{9}}
\,{f_{3}}\,{f_{2}}^{2} + 4112640\,{f_{8}}\,{f_{4}}\,{f_{2}}^{2}
 + 16181760\,{f_{7}}\,{f_{5}}\,{f_{2}}^{2} \\
\mbox{} + 33078528\,{f_{6}}\,{f_{3}}^{2}\,{f_{2}}^{2} + 1650432\,
{f_{6}}\,{f_{4}}\,{f_{2}}^{3} - 748608\,{f_{4}}\,{f_{3}}^{2}\,{f
_{2}}^{3} - 19176640\,{f_{5}}\,{f_{3}}\,{f_{2}}^{4} \\
\mbox{} + 9069312\,{f_{6}}\,{f_{4}}\,{f_{3}}^{2} - 1007104\,{f_{6
}}\,{f_{4}}^{2}\,{f_{2}} + 4305\,{f_{4}}^{4} - 1007\,{f_{2}}^{8}
 + 189665280\,{f_{8}}^{2}=0
\end{array}
$$
$$
\begin{array}{l}
1693440\,{f_{5}}\,{f_{3}}^{2}\,{f_{4}}\,{f_{2}} + 188160\,{f_{8}}
\,{f_{3}}\,{f_{4}}\,{f_{2}} + 23802240\,{f_{6}}\,{f_{5}}\,{f_{4}}
\,{f_{2}} + 412048\,{f_{6}}\,{f_{3}}\,{f_{2}}^{2}\,{f_{4}} \\
\mbox{} - 1011548160\,{f_{6}}^{2}\,{f_{5}} + 284497920\,{f_{9}}\,
{f_{8}} - 167160\,{f_{5}}\,{f_{2}}^{6} + 88200\,{f_{5}}\,{f_{4}}
^{3} + 70560\,{f_{9}}\,{f_{2}}^{4} \\
\mbox{} + 4913121024000\,{f_{5}}^{3}\,{f_{2}} + 85680\,{f_{4}}^{2
}\,{f_{3}}^{3} - 1485711360\,{f_{5}}\,{f_{3}}^{4} + 2184\,{f_{7}}
\,{f_{2}}^{5} \\
\mbox{} - 98262420480\,{f_{7}}\,{f_{5}}^{2} + 526848\,{f_{10}}\,{
f_{7}} + 123312\,{f_{3}}^{3}\,{f_{2}}^{4} + 10536960\,{f_{3}}^{5}
\,{f_{2}} \\
\mbox{} - 31610880\,{f_{8}}\,{f_{3}}^{3} - 2071\,{f_{3}}\,{f_{2}}
^{7} - 493920\,{f_{9}}\,{f_{4}}^{2} + 1818880\,{f_{5}}\,{f_{3}}^{
2}\,{f_{2}}^{3} + 4757\,{f_{3}}\,{f_{2}}^{5}\,{f_{4}} \\
\mbox{} + 74889467520\,{f_{5}}^{2}\,{f_{3}}\,{f_{2}}^{2} - 4368\,
{f_{7}}\,{f_{2}}^{3}\,{f_{4}} - 26342400\,{f_{10}}\,{f_{5}}\,{f_{
2}} + 56448\,{f_{7}}\,{f_{4}}\,{f_{3}}^{2} \\
\mbox{} + 615\,{f_{4}}^{3}\,{f_{3}}\,{f_{2}} + 2184\,{f_{7}}\,{f
_{4}}^{2}\,{f_{2}} - 115080\,{f_{6}}\,{f_{4}}^{2}\,{f_{3}} + 
272832\,{f_{7}}\,{f_{2}}^{2}\,{f_{6}} - 272832\,{f_{7}}\,{f_{6}}
\,{f_{4}} \\
\mbox{} - 208992\,{f_{3}}^{3}\,{f_{2}}^{2}\,{f_{4}} + 2402426880
\,{f_{6}}\,{f_{5}}\,{f_{3}}^{2} - 188160\,{f_{8}}\,{f_{3}}\,{f_{2
}}^{3} + 47040\,{f_{10}}\,{f_{4}}\,{f_{3}} \\
\end{array}
$$

$$
\begin{array}{l}
- 609623613760\,{f_{6}}\,{f_{5}}\,{f_{3}}\,{f_{2}}^{2} - 
14374530240\,{f_{7}}\,{f_{5}}\,{f_{4}}\,{f_{2}} - 19882914240\,{f
_{5}}\,{f_{4}}^{2}\,{f_{3}}\,{f_{2}} \\
\mbox{} + 738439705920\,{f_{6}}\,{f_{5}}\,{f_{4}}\,{f_{3}} + 
302779545600\,{f_{8}}\,{f_{5}}\,{f_{3}}\,{f_{2}} + 17615939040\,{
f_{9}}\,{f_{4}}\,{f_{3}}\,{f_{2}} \\
\mbox{} + 8715924784\,{f_{6}}\,{f_{4}}\,{f_{3}}^{2}\,{f_{2}} + 
26684913920\,{f_{5}}\,{f_{4}}\,{f_{3}}\,{f_{2}}^{3} + 1053441\,{f
_{4}}^{4}\,{f_{2}} - 526187\,{f_{4}}\,{f_{2}}^{7} \\
\mbox{} - 63819497\,{f_{6}}\,{f_{2}}^{6} - 2484310976\,{f_{6}}^{2
}\,{f_{2}}^{3} - 100742040\,{f_{10}}\,{f_{2}}^{4} + 2105815\,{f_{
4}}^{2}\,{f_{2}}^{5} \\
\mbox{} + 18064191090960\,{f_{5}}^{2}\,{f_{2}}^{4} - 2633069\,{f
_{4}}^{3}\,{f_{2}}^{3} + 19792514\,{f_{3}}^{2}\,{f_{2}}^{6} \\
\mbox{} + 2158397472\,{f_{3}}^{4}\,{f_{2}}^{3} - 281048376\,{f_{8
}}\,{f_{2}}^{5} - 807740962560\,{f_{6}}^{2}\,{f_{3}}^{2} + 
29521212\,{f_{4}}^{3}\,{f_{3}}^{2} \\
\mbox{} - 3713751552\,{f_{10}}\,{f_{8}} + 752034689280\,{f_{9}}\,
{f_{3}}^{3} + 724181552640\,{f_{6}}\,{f_{3}}^{4} \\
\mbox{} + 48553512\,{f_{10}}\,{f_{4}}^{2} - 8945581947120\,{f_{5}
}^{2}\,{f_{4}}^{2} - 58978983\,{f_{6}}\,{f_{4}}^{3} \\
\end{array}
$$
$$
\begin{array}{l}
 - 33148029168672082560\,{f_{9}}\,{f_{3}}^{3}\,{f_{2}} - 
1615119956572416000\,{f_{10}}\,{f_{5}}\,{f_{3}}\,{f_{2}} \\
\mbox{} + 159410369241955440\,{f_{9}}\,{f_{4}}\,{f_{3}}\,{f_{2}}
^{2} + 6616249713784757760\,{f_{5}}\,{f_{4}}\,{f_{3}}^{3}\,{f_{2}
} \\
\mbox{} + 74910752807005248\,{f_{7}}\,{f_{6}}\,{f_{3}}\,{f_{2}}^{
2} + 183429726867238152\,{f_{6}}\,{f_{4}}\,{f_{3}}^{2}\,{f_{2}}^{
2} \\
\mbox{} + 243683386620972480\,{f_{8}}\,{f_{4}}\,{f_{3}}^{2}\,{f_{
2}} + 156955879038822080\,{f_{5}}\,{f_{4}}\,{f_{3}}\,{f_{2}}^{4}
 \\
\mbox{} + 33282997663046261760\,{f_{9}}\,{f_{6}}\,{f_{3}}\,{f_{2}
} - 161855888224392040\,{f_{5}}\,{f_{4}}^{2}\,{f_{3}}\,{f_{2}}^{2
} \\
\mbox{} + 2160282922459795840\,{f_{6}}\,{f_{5}}\,{f_{3}}\,{f_{2}}
^{3} - 74910752807005248\,{f_{7}}\,{f_{6}}\,{f_{4}}\,{f_{3}} \\
\mbox{} + 36993024884145012480\,{f_{8}}\,{f_{5}}\,{f_{3}}\,{f_{2}
}^{2} + 599655040942776\,{f_{7}}\,{f_{4}}^{2}\,{f_{3}}\,{f_{2}}
 \\
\mbox{} - 1199310081885552\,{f_{7}}\,{f_{4}}\,{f_{3}}\,{f_{2}}^{3
} - 14685726803748000\,{f_{7}}\,{f_{5}}\,{f_{4}}\,{f_{2}}^{2} \\
\end{array}
$$

$$
\begin{array}{l}
- 376564832\,{f_{10}}\,{f_{3}}\,{f_{2}}^{3} + 69222949461120\,{f
_{5}}^{2}\,{f_{3}}\,{f_{2}}^{3} - 77762764800\,{f_{8}}\,{f_{5}}\,
{f_{3}}^{2} \\
\mbox{} + 4793428080\,{f_{9}}\,{f_{4}}\,{f_{3}}^{2} + 1131285120
\,{f_{6}}\,{f_{4}}\,{f_{3}}^{3} + 11070393600\,{f_{7}}\,{f_{6}}\,
{f_{3}}^{2} \\
\mbox{} - 767151000\,{f_{5}}\,{f_{4}}^{2}\,{f_{3}}^{2} - 
2118060000\,{f_{6}}\,{f_{5}}\,{f_{4}}^{2} - 156548712\,{f_{3}}^{3
}\,{f_{2}}^{3}\,{f_{4}} \\
\mbox{} - 411521160\,{f_{8}}\,{f_{3}}\,{f_{2}}^{4} + 395371200\,{
f_{8}}\,{f_{7}}\,{f_{2}}^{2} + 267446865\,{f_{5}}\,{f_{4}}\,{f_{2
}}^{5} \\
\mbox{} - 186747045\,{f_{5}}\,{f_{4}}^{2}\,{f_{2}}^{3} - 
26327253120\,{f_{10}}\,{f_{5}}\,{f_{2}}^{2} + 217319\,{f_{7}}\,{f
_{4}}^{2}\,{f_{2}}^{2} \\
\mbox{} - 257277368\,{f_{6}}\,{f_{3}}\,{f_{2}}^{5} - 20774130720
\,{f_{6}}\,{f_{5}}\,{f_{2}}^{4} - 2325763\,{f_{7}}\,{f_{2}}^{4}\,
{f_{4}} \\
\mbox{} + 263393312\,{f_{7}}\,{f_{2}}^{3}\,{f_{6}} - 54495168\,{f
_{7}}\,{f_{3}}^{2}\,{f_{2}}^{3} - 7178672200\,{f_{5}}\,{f_{3}}^{2
}\,{f_{2}}^{4} \\
\mbox{} - 21424717440\,{f_{3}}^{3}\,{f_{6}}\,{f_{2}}^{2} + 
285866280\,{f_{9}}\,{f_{2}}^{3}\,{f_{4}} - 7170770880\,{f_{8}}\,{
f_{2}}^{3}\,{f_{5}} \\
\end{array}
$$

\newpage
\section*{Appendix}
The generating invariants $f_4,f_8,f_{12},f_{18},$ in the case $n=5.$ Here we denote the variable $x_0$ by $t$.
$$
\begin{array}{l}
{f_{4}} := t^{2}\,{u_{5}}^{2} + 4\,t\,{u_{2}}\,{u_{3}}\,{u_{5}}
 + 48\,{u_{4}}\,{u_{2}}^{3} - 32\,{u_{2}}^{2}\,{u_{3}}^{2} + 16\,
t\,{u_{2}}\,{u_{4}}^{2} - 12\,t\,{u_{3}}^{2}\,{u_{4}},\mbox{\phantom{************} }\\
\mbox{\phantom{******}}\\
{f_{8}} := 22\,t^{2}\,{u_{2}}^{2}\,{u_{4}}^{4} + 8\,{u_{2}}^{4}\,
{u_{3}}^{4} + 18\,{u_{2}}^{6}\,{u_{4}}^{2} - 24\,{u_{2}}^{5}\,{u
_{3}}^{2}\,{u_{4}} - 38\,t\,{u_{2}}^{4}\,{u_{4}}^{3} + 18\,t^{2}
\,{u_{3}}^{4}\,{u_{4}}^{2}- \\
\mbox{} - 27\,t^{2}\,{u_{5}}\,{u_{3}}^{5} - 2\,t^{3}\,{u_{4}}^{5}
 + 93\,t\,{u_{2}}^{4}\,{u_{3}}\,{u_{4}}\,{u_{5}} - 34\,t^{2}\,{u
_{2}}^{2}\,{u_{3}}\,{u_{5}}\,{u_{4}}^{2} + 12\,t^{2}\,{u_{2}}^{3}
\,{u_{5}}^{2}\,{u_{4}}+ \\
\mbox{} + 78\,t^{2}\,{u_{2}}\,{u_{3}}^{3}\,{u_{5}}\,{u_{4}} - 27
\,t\,{u_{5}}^{2}\,{u_{2}}^{5} - 21\,t^{2}\,{u_{2}}^{2}\,{u_{5}}^{
2}\,{u_{3}}^{2} - 48\,t^{2}\,{u_{2}}\,{u_{3}}^{2}\,{u_{4}}^{3}+8\,t\,{u_{2}}^{3}\,{u_{3}}^{2}\,{u_{4}}^{2}
 \\
\mbox{}  - t^{3}\,{u
_{2}}\,{u_{5}}^{2}\,{u_{4}}^{2} + 6\,t\,{u_{2}}^{2}\,{u_{3}}^{4}
\,{u_{4}} + t^{3}\,{u_{2}}\,{u_{3}}\,{u_{5}}^{3} - 42\,t\,{u_{2}}
^{3}\,{u_{5}}\,{u_{3}}^{3} \mbox{} - 3\,t^{3}\,{u_{5}}^{2}\,{u_{3}}^{2}\,{u_{4}} + 5\,t^{3}
\,{u_{3}}\,{u_{5}}\,{u_{4}}^{3}\\
\mbox{\phantom{******}}\\
{f_{12}} := 252\,t^{2}\,{u_{3}}^{2}\,{u_{5}}^{2}\,{u_{4}}\,{u_{2}
}^{5} - 2010\,t\,{u_{2}}^{6}\,{u_{5}}\,{u_{3}}^{3}\,{u_{4}} - 470
\,t^{2}\,{u_{2}}^{4}\,{u_{5}}\,{u_{3}}^{3}\,{u_{4}}^{2} \\
\end{array}
$$
\newpage
$$
\begin{array}{l}
\mbox{} + 1890\,t^{2}\,{u_{5}}^{3}\,{u_{2}}^{11}\,{u_{4}}^{2} - 
25\,t^{4}\,{u_{2}}^{7}\,{u_{5}}^{3}\,{u_{4}}^{4} + 114960\,t^{2}
\,{u_{2}}^{7}\,{u_{3}}^{5}\,{u_{4}}^{4} + 2481\,t^{2}\,{u_{2}}^{5
}\,{u_{3}}^{10}\,{u_{5}} \\
\mbox{} - 10\,t^{7}\,{u_{5}}^{4}\,{u_{3}}^{3}\,{u_{4}}^{4} + 40\,
t^{6}\,{u_{5}}^{3}\,{u_{4}}^{6}\,{u_{2}}^{3} + 60\,t^{4}\,{u_{3}}
^{2}\,{u_{5}}^{5}\,{u_{2}}^{7} + 10\,t^{6}\,{u_{5}}^{5}\,{u_{4}}
^{3}\,{u_{2}}^{4}+ \\
\mbox{} - 3900\,t^{4}\,{u_{3}}^{5}\,{u_{4}}^{3}\,{u_{5}}^{2}\,{u
_{2}}^{4} - 265\,t^{6}\,{u_{2}}\,{u_{3}}^{4}\,{u_{5}}^{3}\,{u_{4}
}^{4} - 1240\,t^{4}\,{u_{3}}^{4}\,{u_{4}}^{5}\,{u_{2}}^{4}\,{u_{5
}} \\
\mbox{} + 10595\,t\,{u_{5}}\,{u_{2}}^{10}\,{u_{3}}^{4}\,{u_{4}}^{
2} + 30\,t^{5}\,{u_{5}}^{5}\,{u_{3}}^{2}\,{u_{2}}^{5}\,{u_{4}} - 
17520\,t\,{u_{5}}\,{u_{2}}^{11}\,{u_{3}}^{2}\,{u_{4}}^{3} \\
\mbox{} + 270\,t\,{u_{5}}^{2}\,{u_{2}}^{12}\,{u_{3}}\,{u_{4}}^{2}
 - 1650\,t\,{u_{2}}^{9}\,{u_{3}}^{6}\,{u_{5}}\,{u_{4}} + 43605\,t
^{2}\,{u_{5}}\,{u_{2}}^{9}\,{u_{3}}^{2}\,{u_{4}}^{4} \\
\mbox{} - 360\,t^{5}\,{u_{2}}\,{u_{3}}^{8}\,{u_{5}}^{3}\,{u_{4}}
 + 1320\,t^{5}\,{u_{2}}\,{u_{3}}^{7}\,{u_{5}}^{2}\,{u_{4}}^{3} - 
1110\,t^{5}\,{u_{2}}\,{u_{4}}^{5}\,{u_{3}}^{6}\,{u_{5}} \\
\mbox{} + 110\,t^{6}\,{u_{2}}\,{u_{3}}^{5}\,{u_{5}}^{4}\,{u_{4}}
^{2} - 48360\,t^{3}\,{u_{2}}^{5}\,{u_{3}}^{6}\,{u_{4}}^{3}\,{u_{5
}} - 7590\,t^{2}\,{u_{5}}^{2}\,{u_{2}}^{10}\,{u_{3}}\,{u_{4}}^{3}
 \\
\mbox{} + 2920\,t^{4}\,{u_{2}}^{3}\,{u_{5}}^{2}\,{u_{3}}^{7}\,{u
_{4}}^{2} - 65\,t^{6}\,{u_{2}}\,{u_{4}}^{8}\,{u_{3}}^{2}\,{u_{5}}
 + 1275\,t^{4}\,{u_{3}}^{2}\,{u_{4}}^{6}\,{u_{2}}^{5}\,{u_{5}}
 \\
\end{array}
$$
\newpage
The generating invariants $f_2,f_4,f_6,f_{10},f_{15}.$ in the case  $n=6.$
$$
\begin{array}{l}
{f_{2}} :=  - t\,{u_{6}} - 15\,{u_{2}}\,{u_{4}} + 10\,{u_{3}}^{2}.\\
\mbox{\phantom{*********}}\\
f_4:=1200\,{u_{2}}\,{u_{3}}^{2}\,{u_{4}} + 300\,t\,{u_{5}}^{2}\,{u_{2}
} + 300\,{u_{6}}\,{u_{2}}^{3} + 320\,t\,{u_{6}}\,{u_{3}}^{2} - 
600\,t\,{u_{3}}\,{u_{5}}\,{u_{4}} -  \\
\mbox{} -330\,t\,{u_{2}}\,{u_{6}}\,{u
_{4}} - 400\,{u_{3}}^{4} - t^{2}\,{u_{6}}^{2} + 300\,t\,{u_{4}}
^{3} - 525\,{u_{2}}^{2}\,{u_{4}}^{2} - 600\,{u_{2}}^{2}\,{u_{3}}
\,{u_{5}} 
\end{array}
$$
$$
\begin{array}{l}
f_6=12\,t^{2}\,{u_{3}}\,{u_{5}}^{3} + 16\,t\,{u_{3}}^{2}\,{u_{4}}^{3}
 + 8\,t\,{u_{2}}^{3}\,{u_{6}}^{2} + 36\,{u_{2}}^{3}\,{u_{4}}^{3}
 - 81\,{u_{5}}^{2}\,{u_{2}}^{4} - 24\,t\,{u_{2}}\,{u_{4}}^{4} \\
\mbox{} - 337\,{u_{2}}^{2}\,{u_{3}}^{2}\,{u_{4}}^{2} + 16\,{u_{6}
}\,{u_{2}}^{3}\,{u_{3}}^{2} + 318\,{u_{2}}^{3}\,{u_{3}}\,{u_{5}}
\,{u_{4}} - 152\,{u_{2}}^{2}\,{u_{3}}^{3}\,{u_{5}} - 24\,{u_{2}}
^{4}\,{u_{6}}\,{u_{4}} \\
\mbox{} - 16\,t\,{u_{6}}\,{u_{3}}^{4} + 288\,{u_{2}}\,{u_{3}}^{4}
\,{u_{4}} - 10\,t^{2}\,{u_{3}}\,{u_{5}}\,{u_{4}}\,{u_{6}} - 46\,t
\,{u_{2}}^{2}\,{u_{3}}\,{u_{5}}\,{u_{6}} + 30\,t\,{u_{2}}^{2}\,{u
_{4}}\,{u_{5}}^{2} \\
\mbox{} - 64\,{u_{3}}^{6} + 8\,t^{2}\,{u_{6}}\,{u_{4}}^{3} - 6\,t
\,{u_{2}}\,{u_{3}}\,{u_{5}}\,{u_{4}}^{2} - t^{2}\,{u_{3}}^{2}\,{u
_{6}}^{2} - 9\,t^{2}\,{u_{5}}^{2}\,{u_{4}}^{2} - 4\,t^{2}\,{u_{5}
}^{2}\,{u_{2}}\,{u_{6}} \\
\mbox{} + 50\,t\,{u_{2}}\,{u_{3}}^{2}\,{u_{6}}\,{u_{4}} + 4\,t\,{
u_{2}}\,{u_{3}}^{2}\,{u_{5}}^{2} + 4\,t^{2}\,{u_{2}}\,{u_{4}}\,{u
_{6}}^{2} - 8\,t\,{u_{3}}^{3}\,{u_{5}}\,{u_{4}} - 8\,t\,{u_{2}}^{
2}\,{u_{6}}\,{u_{4}}^{2}.\\
\phantom{**************}\\
f_{10}=- 960\,t\,{u_{2}}^{2}\,{u_{4}}^{5}\,{u_{3}}^{2} - 126\,t\,{u_{4}
}^{4}\,{u_{2}}^{4}\,{u_{6}} + 900\,t\,{u_{2}}\,{u_{4}}^{4}\,{u_{3
}}^{4} - 510\,t^{2}\,{u_{3}}^{4}\,{u_{5}}^{2}\,{u_{4}}^{2} + \\
\mbox{} + 12\,
t^{2}\,{u_{4}}^{3}\,{u_{3}}^{4}\,{u_{6}}+ 243\,t\,{u_{2}}^{5}\,{u_{5}}^{4} - 84\,t^{2}\,{u_{3}}^{
5}\,{u_{5}}\,{u_{4}}\,{u_{6}} - 96\,t\,{u_{2}}^{4}\,{u_{6}}^{2}\,
{u_{3}}^{2}\,{u_{4}} - \\
\end{array}
$$
\newpage
$$
\begin{array}{l}
f_{15}=432\,t^{4}\,{u_{4}}^{10}\,{u_{5}} - 2816\,t^{4}\,{u_{5}}^{6}\,{u
_{3}}^{5} - 8000\,{u_{2}}^{6}\,{u_{3}}^{7}\,{u_{6}}^{2} + 8000\,t
^{2}\,{u_{4}}^{6}\,{u_{3}}^{7} - 18225\,{u_{2}}^{9}\,{u_{4}}^{3}
\,{u_{5}}^{3} \\
\mbox{} - 756\,t^{3}\,{u_{5}}^{7}\,{u_{2}}^{5} - 12800\,{u_{2}}^{
6}\,{u_{3}}^{6}\,{u_{5}}^{3} - 19683\,{u_{2}}^{10}\,{u_{5}}^{5}
 + 1600\,{u_{2}}^{9}\,{u_{6}}^{3}\,{u_{3}}^{3} \\
\mbox{} - 1600\,t^{3}\,{u_{4}}^{9}\,{u_{3}}^{3} + 36\,t^{5}\,{u_{
4}}^{6}\,{u_{5}}^{3}\,{u_{6}} - 31860\,{u_{2}}^{8}\,{u_{3}}^{3}\,
{u_{5}}^{4} - 32768\,t^{2}\,{u_{3}}^{10}\,{u_{5}}^{3} \\
\mbox{} + 9\,t^{5}\,{u_{4}}^{5}\,{u_{5}}^{5} + 2160\,{u_{2}}^{11}
\,{u_{6}}^{3}\,{u_{5}} - 5670\,t^{4}\,{u_{2}}^{2}\,{u_{3}}\,{u_{4
}}^{4}\,{u_{6}}^{2}\,{u_{5}}^{2} - 2960\,t\,{u_{2}}^{6}\,{u_{6}}
^{3}\,{u_{3}}^{5} \\
\mbox{} + 105300\,t^{2}\,{u_{2}}^{5}\,{u_{4}}\,{u_{3}}^{2}\,{u_{5
}}^{5} + 675\,t^{3}\,{u_{2}}^{4}\,{u_{4}}^{3}\,{u_{5}}^{5} + 2960
\,t^{3}\,{u_{4}}^{6}\,{u_{3}}^{5}\,{u_{6}} \\
\mbox{} - 132\,t^{5}\,{u_{2}}\,{u_{4}}^{5}\,{u_{6}}^{3}\,{u_{5}}
 - 27000\,t^{2}\,{u_{2}}^{3}\,{u_{4}}^{9}\,{u_{3}} - 6960\,t^{4}
\,{u_{2}}\,{u_{3}}^{3}\,{u_{4}}^{3}\,{u_{6}}^{2}\,{u_{5}}^{2} \\
\mbox{} + 15\,t^{4}\,{u_{5}}\,{u_{2}}^{5}\,{u_{6}}^{5} - 945\,t^{
3}\,{u_{2}}^{3}\,{u_{4}}^{6}\,{u_{5}}^{3} - 45\,t^{4}\,{u_{2}}^{3
}\,{u_{5}}^{7}\,{u_{4}} - 50625\,{u_{2}}^{8}\,{u_{4}}^{5}\,{u_{6}
}\,{u_{3}} \\
\mbox{} + 22536\,t^{2}\,{u_{4}}\,{u_{3}}^{3}\,{u_{5}}^{2}\,{u_{2}
}^{5}\,{u_{6}}^{2} - 3\,t^{5}\,{u_{2}}^{2}\,{u_{5}}^{7}\,{u_{6}}
 + 464\,t^{4}\,{u_{4}}^{6}\,{u_{6}}^{2}\,{u_{3}}^{3} \\
\end{array}
$$
\newpage 
$$
\begin{array}{l}
\mbox{} + 16186\,t^{3}\,{u_{2}}^{3}\,{u_{4}}^{6}\,{u_{6}}^{2}\,{u
_{3}} + 336\,t^{5}\,{u_{4}}^{2}\,{u_{3}}^{2}\,{u_{6}}\,{u_{5}}^{5
} + 8640\,t^{4}\,{u_{4}}^{2}\,{u_{3}}^{4}\,{u_{5}}^{5} \\
\mbox{} - 30375\,t\,{u_{2}}^{6}\,{u_{4}}^{7}\,{u_{5}} + t^{6}\,{u
_{5}}^{6}\,{u_{3}}\,{u_{6}}^{2} - 615\,t^{3}\,{u_{2}}^{4}\,{u_{3}
}\,{u_{5}}^{4}\,{u_{4}}^{2}\,{u_{6}} + 30375\,{u_{2}}^{8}\,{u_{4}
}^{4}\,{u_{5}}^{2}\,{u_{3}} \\
\mbox{} - 48600\,t\,{u_{2}}^{6}\,{u_{6}}\,{u_{3}}^{2}\,{u_{4}}^{4
}\,{u_{5}} - 1800\,t^{4}\,{u_{2}}^{2}\,{u_{3}}^{3}\,{u_{4}}\,{u_{
6}}^{3}\,{u_{5}}^{2} + 69\,t^{5}\,{u_{2}}\,{u_{3}}\,{u_{4}}^{4}\,
{u_{6}}^{4} \\
\mbox{} - 61560\,t\,{u_{2}}^{7}\,{u_{3}}^{2}\,{u_{5}}^{5} - 13932
\,t\,{u_{2}}^{9}\,{u_{6}}^{2}\,{u_{5}}^{3} - 42000\,t^{2}\,{u_{2}
}\,{u_{4}}^{7}\,{u_{3}}^{5} \\
\mbox{} - 20400\,{u_{2}}^{8}\,{u_{3}}^{4}\,{u_{6}}^{2}\,{u_{5}}
 - 42525\,t\,{u_{2}}^{7}\,{u_{4}}^{4}\,{u_{5}}^{3} - 11\,t^{6}\,{
u_{4}}^{3}\,{u_{6}}^{3}\,{u_{5}}^{3} - 64\,t^{5}\,{u_{3}}^{4}\,{u
_{6}}^{3}\,{u_{5}}^{3} \\
\mbox{} + 43740\,{u_{2}}^{10}\,{u_{6}}\,{u_{5}}^{3}\,{u_{4}} + 
6000\,{u_{2}}^{6}\,{u_{4}}^{2}\,{u_{3}}^{5}\,{u_{5}}^{2} - 960\,t
^{4}\,{u_{3}}^{2}\,{u_{2}}^{2}\,{u_{5}}^{7} \\
\mbox{} + 61200\,{u_{3}}^{4}\,{u_{4}}\,{u_{5}}^{3}\,{u_{2}}^{7}
 - 3600\,{u_{2}}^{10}\,{u_{6}}^{3}\,{u_{3}}\,{u_{4}} + 42000\,{u
_{2}}^{7}\,{u_{6}}^{2}\,{u_{3}}^{5}\,{u_{4}} \\
\mbox{} - 12960\,{u_{2}}^{10}\,{u_{6}}^{2}\,{u_{3}}\,{u_{5}}^{2}
 - 2400\,{u_{2}}^{7}\,{u_{6}}\,{u_{5}}^{2}\,{u_{3}}^{5} + 27000\,
{u_{2}}^{9}\,{u_{4}}^{3}\,{u_{6}}^{2}\,{u_{3}} \\
\end{array}
$$
\newpage 
$$
\begin{array}{l}
\mbox{} - 91600\,t\,{u_{4}}^{3}\,{u_{6}}^{2}\,{u_{2}}^{6}\,{u_{3}
}^{3} + 45600\,t^{3}\,{u_{2}}^{2}\,{u_{3}}^{4}\,{u_{5}}\,{u_{4}}
^{3}\,{u_{6}}^{2} + 1665\,t^{4}\,{u_{2}}^{2}\,{u_{3}}\,{u_{4}}^{2
}\,{u_{5}}^{6} \\
\mbox{} + 18\,t^{5}\,{u_{4}}\,{u_{2}}^{2}\,{u_{6}}^{2}\,{u_{5}}^{
5} + 230850\,t\,{u_{2}}^{6}\,{u_{4}}^{5}\,{u_{3}}\,{u_{5}}^{2} + 
482400\,t^{2}\,{u_{2}}^{2}\,{u_{3}}^{5}\,{u_{4}}^{4}\,{u_{5}}^{2}
 \\
\mbox{} + 184320\,t\,{u_{2}}^{4}\,{u_{3}}^{6}\,{u_{5}}^{3}\,{u_{4
}} - 18912\,t^{3}\,{u_{2}}^{2}\,{u_{3}}^{3}\,{u_{4}}^{5}\,{u_{6}}
^{2} - 6960\,t\,{u_{4}}^{2}\,{u_{3}}^{4}\,{u_{2}}^{5}\,{u_{5}}^{3
} \\
\mbox{} - 9600\,t\,{u_{3}}^{6}\,{u_{4}}^{2}\,{u_{2}}^{4}\,{u_{6}}
\,{u_{5}} + 12960\,t\,{u_{2}}^{7}\,{u_{6}}^{3}\,{u_{3}}^{3}\,{u_{
4}} - 243000\,t\,{u_{2}}^{5}\,{u_{4}}^{6}\,{u_{3}}^{2}\,{u_{5}}
 \\
\mbox{} - 201\,t^{4}\,{u_{5}}^{2}\,{u_{6}}^{4}\,{u_{3}}\,{u_{2}}
^{4} - 155520\,t\,{u_{2}}^{5}\,{u_{6}}\,{u_{5}}^{2}\,{u_{3}}^{5}
\,{u_{4}} + 222000\,t\,{u_{4}}^{5}\,{u_{3}}^{4}\,{u_{2}}^{4}\,{u
_{5}} \\
\mbox{} + 190920\,t\,{u_{2}}^{6}\,{u_{4}}^{2}\,{u_{5}}^{2}\,{u_{3
}}^{3}\,{u_{6}} + 31500\,t\,{u_{2}}^{7}\,{u_{6}}^{2}\,{u_{3}}^{2}
\,{u_{5}}\,{u_{4}}^{2} + 52800\,t\,{u_{4}}^{2}\,{u_{3}}^{5}\,{u_{
2}}^{5}\,{u_{6}}^{2} \\
\mbox{} + 213120\,t\,{u_{4}}\,{u_{2}}^{6}\,{u_{5}}^{4}\,{u_{3}}^{
3} - 9600\,t\,{u_{2}}^{4}\,{u_{3}}^{7}\,{u_{6}}^{2}\,{u_{4}} - 
39900\,t^{3}\,{u_{3}}^{2}\,{u_{2}}^{3}\,{u_{4}}^{4}\,{u_{6}}^{2}
\,{u_{5}}.
\end{array}
$$
\newpage
Generating set $f_2,f_3,f_4,f_{5},f_{6},f_7,f_8,f_9,f_{10}$ for the invariant ring  $R_8^{SL_2}.$
$$
\begin{array}{l}
f_{2}= - 56\,{u_{3}}\,{u_{5}} + 28\,{u_{2}}\,{u_{6}} + t\,{u_{8
}} + 35\,{u_{4}}^{2}.\\
\mbox{\phantom{*****}}\\
f_{3}=t\,{u_{4}}\,{u_{8}} - 8\,{u_{2}}\,{u_{3}}\,{u_{7}} - 36\,
{u_{3}}\,{u_{4}}\,{u_{5}} - 22\,{u_{2}}\,{u_{4}}\,{u_{6}} - 4\,t
\,{u_{5}}\,{u_{7}} + 15\,{u_{4}}^{3} + 3\,t\,{u_{6}}^{2} + 3\,{u
_{2}}^{2}\,{u_{8}} \\
\mbox{} + 24\,{u_{2}}\,{u_{5}}^{2} + 24\,{u_{3}}^{2}\,{u_{6}}.\\
\mbox{\phantom{*****}}\\
{f_{4}}=168\,{u_{2}}^{2}\,{u_{4}}\,{u_{8}} + 2408\,{u_{3}}^{2}\,{
u_{5}}^{2} + 504\,{u_{2}}^{2}\,{u_{5}}\,{u_{7}} - 1176\,{u_{3}}^{
2}\,{u_{4}}\,{u_{6}} + 168\,t\,{u_{4}}\,{u_{6}}^{2} \\
\mbox{} - 112\,{u_{2}}\,{u_{3}}^{2}\,{u_{8}} + 56\,t\,{u_{2}}\,{u
_{7}}^{2} - 1288\,{u_{2}}\,{u_{3}}\,{u_{4}}\,{u_{7}} + 56\,t\,{u
_{3}}\,{u_{5}}\,{u_{8}} + 56\,t\,{u_{4}}\,{u_{5}}\,{u_{7}} \\
\mbox{} - 168\,t\,{u_{3}}\,{u_{6}}\,{u_{7}} + 665\,{u_{4}}^{4} + 
672\,{u_{3}}^{3}\,{u_{7}} - 1624\,{u_{2}}\,{u_{3}}\,{u_{5}}\,{u_{
6}} - 42\,t\,{u_{4}}^{2}\,{u_{8}} \\
\end{array}
$$
\newpage
$$
\begin{array}{l}
{f_{7}}=29400\,{u_{6}}^{2}\,{u_{3}}^{4}\,{u_{4}} - 75\,{u_{8}}^{2
}\,{u_{2}}^{4}\,{u_{4}} + 31800\,{u_{5}}^{4}\,{u_{2}}\,{u_{3}}^{2
} - 75\,{u_{6}}^{4}\,t^{2}\,{u_{4}} + 50\,{u_{8}}^{2}\,{u_{2}}^{3
}\,{u_{3}}^{2} \\
\mbox{} + 6000\,{u_{5}}^{5}\,t\,{u_{3}} + 5600\,{u_{7}}^{2}\,{u_{
3}}^{4}\,{u_{2}} - 6875\,{u_{6}}^{2}\,t\,{u_{4}}^{4} + 3528\,{u_{
5}}^{3}\,{u_{7}}\,{u_{2}}^{3} \\
\mbox{} + 31800\,{u_{5}}^{2}\,{u_{6}}\,{u_{3}}^{4} + 3212\,{u_{6}
}^{4}\,t\,{u_{2}}^{2} - 63000\,{u_{5}}\,{u_{4}}^{5}\,{u_{3}} + 
3000\,{u_{8}}\,{u_{4}}^{5}\,t + 50\,{u_{5}}^{2}\,{u_{6}}^{3}\,t^{
2} \\
\mbox{} + 6250\,{u_{7}}^{2}\,{u_{2}}^{3}\,{u_{4}}^{2} - 1400\,{u
_{6}}\,{u_{7}}^{2}\,{u_{2}}^{4} - 6875\,{u_{8}}\,{u_{2}}^{2}\,{u
_{4}}^{4} + 2\,{u_{8}}^{3}\,t^{2}\,{u_{2}}^{2} \\
\mbox{} + 31125\,{u_{5}}^{2}\,{u_{2}}\,{u_{4}}^{4} - 12632\,{u_{6
}}^{3}\,{u_{2}}^{3}\,{u_{4}} + 10792\,{u_{5}}^{2}\,{u_{6}}^{2}\,{
u_{2}}^{3} - 20375\,{u_{6}}\,{u_{4}}^{5}\,{u_{2}} \\
\mbox{} - 236\,{u_{5}}\,{u_{7}}\,{u_{8}}\,t\,{u_{2}}\,{u_{3}}^{2}
 + 150\,{u_{8}}^{2}\,t^{2}\,{u_{4}}^{3} - 16800\,{u_{6}}\,{u_{7}}
\,{u_{3}}^{5} - 3750\,{u_{5}}^{4}\,t\,{u_{4}}^{2} \\
\mbox{} + 2\,{u_{6}}^{2}\,{u_{8}}^{2}\,t^{3} + 29400\,{u_{5}}^{4}
\,{u_{2}}^{2}\,{u_{4}} + 10792\,{u_{6}}^{3}\,{u_{2}}^{2}\,{u_{3}}
^{2} + 107325\,{u_{5}}^{2}\,{u_{3}}^{2}\,{u_{4}}^{3} \\
\mbox{} + 3212\,{u_{6}}^{2}\,{u_{8}}\,{u_{2}}^{4} + 31125\,{u_{6}
}\,{u_{3}}^{2}\,{u_{4}}^{4} - 4072\,{u_{6}}\,{u_{7}}^{2}\,t\,{u_{
2}}\,{u_{3}}^{2} - 3750\,{u_{8}}\,{u_{4}}^{2}\,{u_{3}}^{4} \\
\end{array}
$$
\newpage 
$$
\begin{array}{l}
{f_{8}}=480\,{u_{6}}\,{u_{4}}^{5}\,{u_{3}}^{2} + 114\,{u_{5}}^{4}
\,{u_{6}}^{2}\,t^{2} - 25\,{u_{5}}^{4}\,t\,{u_{4}}^{3} - 1056\,{u
_{5}}^{2}\,{u_{3}}^{2}\,{u_{4}}^{4} - {u_{7}}^{4}\,t^{3}\,{u_{4}}
 - {u_{6}}^{3}\,{u_{7}}^{2}\,t^{3} \\
\mbox{} + 3741\,{u_{6}}^{2}\,{u_{2}}^{2}\,{u_{4}}^{4} + 114\,{u_{
8}}^{2}\,{u_{3}}^{4}\,{u_{2}}^{2} + 480\,{u_{5}}^{2}\,{u_{4}}^{5}
\,{u_{2}} + 74\,{u_{6}}^{2}\,{u_{7}}\,{u_{8}}\,t^{2}\,{u_{2}}\,{u
_{3}} \\
\mbox{} + 1368\,{u_{5}}^{3}\,{u_{3}}^{3}\,{u_{4}}^{2} + 1536\,{u
_{5}}\,{u_{6}}^{2}\,{u_{3}}^{5} - 888\,{u_{5}}^{2}\,{u_{7}}\,{u_{
3}}^{5} - 40\,{u_{5}}^{6}\,t\,{u_{2}} - 40\,{u_{6}}\,{u_{8}}\,{u
_{3}}^{6} \\
\mbox{} + 174\,{u_{8}}^{2}\,{u_{2}}^{4}\,{u_{4}}^{2} - 940\,{u_{7
}}\,{u_{4}}^{4}\,{u_{3}}^{3} - 365\,{u_{8}}\,{u_{4}}^{5}\,{u_{2}}
^{2} + 1371\,{u_{6}}^{2}\,{u_{3}}^{4}\,{u_{4}}^{2} + 180\,{u_{8}}
\,{u_{4}}^{6}\,t \\
\mbox{} + {u_{6}}^{4}\,{u_{8}}\,t^{3} + 45\,{u_{7}}^{2}\,{u_{8}}
\,{u_{2}}^{5} + 1521\,{u_{5}}^{2}\,{u_{7}}^{2}\,{u_{2}}^{4} + 4\,
{u_{5}}\,{u_{8}}^{2}\,{u_{7}}\,t^{2}\,{u_{2}}^{2} \\
\mbox{} + 1076\,{u_{6}}^{3}\,{u_{5}}\,t\,{u_{2}}\,{u_{3}}\,{u_{4}
} - 66\,{u_{6}}^{5}\,t^{2}\,{u_{2}} + 1536\,{u_{5}}^{5}\,{u_{2}}
^{2}\,{u_{3}} + 320\,{u_{5}}\,{u_{4}}^{6}\,{u_{3}} \\
\mbox{} + 42\,{u_{5}}^{3}\,{u_{7}}\,{u_{8}}\,t^{2}\,{u_{2}} + 174
\,{u_{6}}^{4}\,t^{2}\,{u_{4}}^{2} - 160\,{u_{5}}^{5}\,{u_{7}}\,t
^{2} + 100\,{u_{7}}^{2}\,{u_{2}}^{3}\,{u_{4}}^{3} \\
\mbox{} - 190\,{u_{5}}\,{u_{6}}\,{u_{8}}\,{u_{3}}^{3}\,{u_{2}}^{2
} + {u_{8}}^{3}\,{u_{2}}^{4}\,t - 3728\,{u_{6}}^{3}\,{u_{2}}^{3}
\,{u_{4}}^{2} - 120\,{u_{7}}^{3}\,{u_{2}}^{4}\,{u_{3}} \\
\end{array}
$$
\newpage 
$$
\begin{array}{l}
\mbox{} - 200\,{u_{6}}\,{u_{7}}\,{u_{8}}\,t\,{u_{2}}^{2}\,{u_{3}}
\,{u_{4}} + 13\,{u_{7}}^{2}\,{u_{8}}\,{u_{2}}^{2}\,{u_{3}}^{2}\,t
 + 2\,{u_{7}}\,{u_{8}}^{2}\,t^{2}\,{u_{2}}\,{u_{4}}\,{u_{3}} \\
\mbox{} + 583\,{u_{6}}\,{u_{7}}^{2}\,{u_{2}}^{2}\,{u_{4}}^{2}\,t
 + 790\,{u_{5}}^{2}\,{u_{6}}\,{u_{7}}\,{u_{3}}^{3}\,t + 3662\,{u
_{5}}^{2}\,{u_{6}}\,{u_{7}}\,{u_{2}}^{3}\,{u_{3}} \\
\mbox{} + 188\,{u_{5}}^{2}\,{u_{6}}\,{u_{7}}\,t\,{u_{2}}\,{u_{3}}
\,{u_{4}} - 602\,{u_{5}}^{3}\,{u_{7}}\,{u_{2}}\,{u_{4}}^{2}\,t - 
5174\,{u_{5}}^{3}\,{u_{6}}\,{u_{2}}^{2}\,{u_{3}}\,{u_{4}} \\
\mbox{} - 198\,{u_{6}}^{3}\,{u_{7}}\,t^{2}\,{u_{3}}\,{u_{4}} - 34
\,{u_{6}}^{3}\,{u_{8}}\,t^{2}\,{u_{2}}\,{u_{4}} + 4138\,{u_{5}}^{
2}\,{u_{6}}^{2}\,{u_{2}}^{3}\,{u_{4}} + 1681\,{u_{5}}^{2}\,{u_{6}
}^{2}\,{u_{2}}^{2}\,{u_{3}}^{2} \\
\mbox{} - 3054\,{u_{5}}^{4}\,{u_{2}}\,{u_{3}}^{2}\,{u_{4}} + 32\,
{u_{6}}^{4}\,{u_{2}}^{2}\,{u_{4}}\,t + 105\,{u_{6}}^{4}\,{u_{2}}
\,{u_{3}}^{2}\,t + 148\,{u_{6}}^{3}\,{u_{8}}\,{u_{2}}^{3}\,t \\
\mbox{} - 210\,{u_{6}}^{3}\,{u_{5}}\,{u_{3}}^{3}\,t - 2904\,{u_{6
}}^{3}\,{u_{5}}\,{u_{2}}^{3}\,{u_{3}} + 1298\,{u_{5}}^{3}\,{u_{6}
}\,{u_{3}}^{3}\,{u_{2}} - 2682\,{u_{5}}^{3}\,{u_{7}}\,{u_{2}}^{3}
\,{u_{4}} \\
\mbox{} - 1594\,{u_{5}}^{3}\,{u_{7}}\,{u_{2}}^{2}\,{u_{3}}^{2} - 
1078\,{u_{5}}^{3}\,{u_{8}}\,{u_{3}}^{3}\,t - 210\,{u_{5}}^{3}\,{u
_{8}}\,{u_{2}}^{3}\,{u_{3}} - 2\,{u_{7}}^{3}\,{u_{5}}\,{u_{2}}^{3
}\,t \\
\mbox{} + 59\,{u_{5}}^{2}\,{u_{8}}^{2}\,t^{2}\,{u_{3}}^{2} - 2760
\,{u_{5}}\,{u_{6}}^{2}\,{u_{7}}\,{u_{2}}^{4} + 105\,{u_{5}}^{2}\,
{u_{6}}\,{u_{8}}\,{u_{2}}^{4} - 244\,{u_{5}}^{2}\,{u_{7}}^{2}\,t
^{2}\,{u_{4}}^{2} \\
\end{array} 
$$

\newpage 
$$
\begin{array}{l}
{f_{9}}=24\,{u_{6}}^{2}\,{u_{7}}^{2}\,{u_{2}}^{5} - 159\,{u_{6}}
^{4}\,{u_{3}}^{4}\,t + 95\,{u_{8}}\,{u_{4}}^{6}\,{u_{2}}^{2} + 10
\,{u_{7}}^{3}\,{u_{5}}\,{u_{2}}^{5} + 136\,{u_{7}}^{3}\,{u_{3}}^{
5}\,t \\
\mbox{} + 66\,{u_{5}}^{6}\,t\,{u_{3}}^{2} + 15\,{u_{7}}^{4}\,{u_{
2}}^{4}\,t - 1670\,{u_{5}}^{3}\,{u_{7}}\,{u_{2}}^{2}\,{u_{3}}^{2}
\,{u_{4}} - 732\,{u_{6}}^{3}\,{u_{2}}^{3}\,{u_{4}}^{3} \\
\mbox{} - 8\,{u_{7}}\,{u_{8}}^{2}\,t^{2}\,{u_{3}}^{3}\,{u_{4}} + 
75\,{u_{5}}^{4}\,t\,{u_{4}}^{4} + 40\,{u_{8}}^{2}\,{u_{2}}^{4}\,{
u_{4}}^{3} + 280\,{u_{6}}^{4}\,{u_{2}}^{3}\,{u_{3}}^{2} + 78\,{u
_{5}}^{2}\,{u_{8}}^{2}\,{u_{2}}^{5} \\
\mbox{} + {u_{7}}^{4}\,t^{3}\,{u_{4}}^{2} - 420\,{u_{5}}^{3}\,{u
_{6}}\,{u_{3}}^{5} + 40\,{u_{6}}^{4}\,t^{2}\,{u_{4}}^{3} - 420\,{
u_{7}}\,{u_{4}}^{5}\,{u_{3}}^{3} + 3456\,{u_{5}}^{2}\,{u_{3}}^{2}
\,{u_{4}}^{5} \\
\mbox{} + 78\,{u_{6}}^{5}\,t^{2}\,{u_{3}}^{2} - 3000\,{u_{5}}^{3}
\,{u_{3}}^{3}\,{u_{4}}^{3} + 72\,{u_{6}}^{3}\,{u_{8}}\,{u_{2}}^{5
} + 75\,{u_{8}}\,{u_{4}}^{4}\,{u_{3}}^{4} + 1200\,{u_{5}}^{4}\,{u
_{2}}^{2}\,{u_{4}}^{3} \\
\mbox{} + {u_{5}}^{4}\,{u_{8}}^{2}\,t^{3} + 240\,{u_{7}}^{2}\,{u
_{3}}^{6}\,{u_{4}} + 66\,{u_{5}}^{2}\,{u_{8}}\,{u_{3}}^{6} + 24\,
{u_{7}}^{3}\,{u_{2}}^{3}\,{u_{3}}^{3} - 50\,{u_{7}}^{2}\,{u_{2}}
^{3}\,{u_{4}}^{4} \\
\mbox{} + 960\,{u_{5}}^{2}\,{u_{4}}^{6}\,{u_{2}} + 1314\,{u_{6}}
^{2}\,{u_{4}}^{5}\,{u_{2}}^{2} + 960\,{u_{6}}\,{u_{4}}^{6}\,{u_{3
}}^{2} + 880\,{u_{6}}\,{u_{8}}\,{u_{2}}^{2}\,{u_{3}}^{2}\,{u_{4}}
^{3} \\
\mbox{} - 1126\,{u_{6}}^{2}\,{u_{8}}\,{u_{2}}\,{u_{3}}^{2}\,{u_{4
}}^{2}\,t - 159\,{u_{5}}^{4}\,{u_{8}}\,{u_{2}}^{4} + 630\,{u_{5}}
^{4}\,{u_{3}}^{4}\,{u_{4}} + 734\,{u_{5}}\,{u_{6}}\,{u_{8}}\,{u_{
2}}^{2}\,{u_{3}}^{3}\,{u_{4}} \\
\end{array}
$$
\newpage
$$
\begin{array}{l}
\mbox{} + 6\,{u_{7}}^{4}\,t^{2}\,{u_{2}}\,{u_{3}}^{2} + 6\,{u_{6}
}^{4}\,t\,{u_{2}}^{2}\,{u_{4}}^{2} + 3900\,{u_{5}}^{4}\,{u_{2}}\,
{u_{3}}^{2}\,{u_{4}}^{2} - 18\,{u_{7}}^{3}\,{u_{6}}\,t^{2}\,{u_{3
}}^{3} \\
\mbox{} - 92\,{u_{5}}^{3}\,{u_{8}}\,{u_{2}}^{2}\,{u_{3}}^{3} - 70
\,{u_{6}}\,{u_{7}}\,t\,{u_{3}}\,{u_{4}}^{5} - 244\,{u_{6}}\,{u_{8
}}\,t\,{u_{2}}\,{u_{4}}^{5} + 540\,{u_{5}}\,{u_{8}}\,{u_{2}}\,{u
_{3}}^{3}\,{u_{4}}^{3} \\
\mbox{} - 240\,{u_{5}}\,{u_{8}}\,{u_{2}}^{2}\,{u_{3}}\,{u_{4}}^{4
} + 550\,{u_{7}}\,{u_{8}}\,{u_{2}}\,{u_{3}}^{5}\,{u_{4}} + 40\,{u
_{7}}\,{u_{8}}\,t\,{u_{3}}^{3}\,{u_{4}}^{3} + 220\,{u_{6}}\,{u_{8
}}\,t\,{u_{3}}^{2}\,{u_{4}}^{4} \\
\mbox{} - 760\,{u_{6}}\,{u_{8}}\,{u_{2}}\,{u_{3}}^{4}\,{u_{4}}^{2
} - 8\,{u_{7}}^{2}\,{u_{8}}\,t^{2}\,{u_{2}}\,{u_{4}}^{3} - 390\,{
u_{5}}\,{u_{6}}\,{u_{7}}\,t\,{u_{2}}\,{u_{4}}^{4} \\
\mbox{} - 2826\,{u_{5}}\,{u_{6}}\,{u_{7}}\,{u_{2}}^{2}\,{u_{3}}^{
2}\,{u_{4}}^{2} + 1510\,{u_{5}}\,{u_{6}}\,{u_{7}}\,{u_{2}}^{3}\,{
u_{4}}^{3} + 398\,{u_{5}}\,{u_{6}}\,{u_{7}}\,t\,{u_{3}}^{2}\,{u_{
4}}^{3} \\
\mbox{} + 618\,{u_{5}}\,{u_{7}}\,{u_{8}}\,{u_{2}}^{3}\,{u_{3}}^{2
}\,{u_{4}} - 770\,{u_{5}}\,{u_{6}}\,{u_{8}}\,t\,{u_{3}}^{3}\,{u_{
4}}^{2} + 22\,{u_{6}}\,{u_{7}}\,{u_{8}}\,t^{2}\,{u_{3}}\,{u_{4}}
^{3} \\
\mbox{} + 346\,{u_{6}}\,{u_{7}}\,{u_{8}}\,{u_{2}}^{4}\,{u_{3}}\,{
u_{4}} + 216\,{u_{5}}\,{u_{7}}\,{u_{8}}\,t\,{u_{2}}^{2}\,{u_{4}}
^{3} + 8\,{u_{7}}\,{u_{8}}^{2}\,t^{2}\,{u_{2}}\,{u_{4}}^{2}\,{u_{
3}} \\
\mbox{} - 310\,{u_{5}}\,{u_{6}}\,{u_{8}}\,{u_{3}}^{5}\,{u_{2}} + 
972\,{u_{5}}\,{u_{6}}\,{u_{8}}\,t\,{u_{2}}\,{u_{4}}^{3}\,{u_{3}}
 - 2\,{u_{8}}^{3}\,t^{2}\,{u_{2}}\,{u_{3}}^{2}\,{u_{4}} \\
\end{array}
$$
\newpage 
$$
\begin{array}{l}
\mbox{} + 28\,{u_{5}}^{3}\,{u_{6}}\,{u_{7}}\,t\,{u_{2}}\,{u_{3}}
^{2} - 54\,{u_{5}}^{3}\,{u_{7}}\,{u_{8}}\,t^{2}\,{u_{2}}\,{u_{4}}
 - 304\,{u_{5}}^{3}\,{u_{6}}\,{u_{8}}\,t\,{u_{2}}^{2}\,{u_{3}}
 \\
\mbox{} - 206\,{u_{5}}^{3}\,{u_{6}}\,{u_{8}}\,t^{2}\,{u_{3}}\,{u
_{4}} + 8\,{u_{5}}^{2}\,{u_{6}}^{2}\,{u_{8}}\,t^{2}\,{u_{2}}\,{u
_{4}} - 24\,{u_{5}}\,{u_{6}}^{2}\,{u_{8}}\,{u_{7}}\,t^{2}\,{u_{2}
}^{2} \\
\mbox{} + 4\,{u_{5}}\,{u_{6}}^{2}\,{u_{8}}\,{u_{7}}\,t^{3}\,{u_{4
}} - 208\,{u_{6}}^{3}\,{u_{5}}\,{u_{7}}\,t^{2}\,{u_{3}}^{2} - 254
\,{u_{6}}^{3}\,{u_{5}}\,{u_{7}}\,{u_{2}}^{3}\,t \\
\mbox{} + 130\,{u_{5}}^{3}\,{u_{6}}\,{u_{7}}\,t^{2}\,{u_{4}}^{2}
 + 104\,{u_{5}}^{2}\,{u_{7}}^{2}\,{u_{6}}\,t^{2}\,{u_{3}}^{2} + 
170\,{u_{5}}^{2}\,{u_{7}}^{2}\,{u_{6}}\,{u_{2}}^{3}\,t - 4\,{u_{7
}}^{3}\,{u_{5}}\,{u_{6}}\,t^{3}\,{u_{4}} \\
\mbox{} - 86\,{u_{5}}^{2}\,{u_{7}}^{2}\,{u_{6}}\,t^{2}\,{u_{2}}\,
{u_{4}} + 102\,{u_{5}}^{2}\,{u_{6}}^{2}\,{u_{8}}\,t^{2}\,{u_{3}}
^{2} + 116\,{u_{5}}^{2}\,{u_{6}}^{2}\,{u_{8}}\,{u_{2}}^{3}\,t \\
\mbox{} - 46\,{u_{6}}^{3}\,{u_{5}}\,{u_{8}}\,t^{2}\,{u_{2}}\,{u_{
3}} - 260\,{u_{7}}^{3}\,{u_{5}}\,t\,{u_{2}}^{3}\,{u_{4}} + 10\,{u
_{5}}^{4}\,{u_{8}}\,{u_{6}}\,t^{2}\,{u_{2}}.\\

\mbox{\phantom{*****}}\\

{f_{10}}=380020\,{u_{8}}\,{u_{4}}^{7}\,{u_{2}}^{2} + 1876355\,{u

_{6}}^{????}\,{u_{4}}^{6}\,{u_{2}}^{2} + 179985\,{u_{8}}^{2}\,{u_{2}
}^{4}\,{u_{4}}^{4} + 2056620\,{u_{5}}^{2}\,{u_{4}}^{6}\,{u_{3}}^{
2} \\
\end{array}
$$
\newpage 
$$
\begin{array}{l}
\mbox{} - 734254\,{u_{5}}^{3}\,{u_{6}}^{2}\,{u_{8}}\,t^{2}\,{u_{2
}}\,{u_{3}} - 3093177\,{u_{5}}^{2}\,{u_{8}}\,{u_{7}}\,{u_{2}}^{4}
\,{u_{3}}\,{u_{4}} \\
\mbox{} + 1313\,{u_{7}}\,{u_{8}}^{2}\,{u_{5}}\,t^{2}\,{u_{2}}^{2}
\,{u_{4}}^{2} + 334\,{u_{7}}^{3}\,{u_{5}}\,{u_{8}}\,{u_{3}}^{2}\,
t^{3} - 1457739\,{u_{6}}^{3}\,{u_{5}}\,{u_{7}}\,{u_{3}}^{2}\,t^{2
}\,{u_{4}} \\
\mbox{} + 4322715\,{u_{5}}^{3}\,{u_{8}}\,{u_{3}}\,t\,{u_{4}}^{3}
\,{u_{2}} - 738540\,{u_{5}}\,{u_{6}}\,{u_{8}}^{2}\,{u_{2}}^{5}\,{
u_{3}} - 703440\,{u_{5}}^{4}\,{u_{8}}\,t\,{u_{2}}^{2}\,{u_{4}}^{2
} \\
\mbox{} + 175\,{u_{5}}\,{u_{6}}^{2}\,{u_{8}}^{2}\,{u_{3}}\,t^{3}
\,{u_{4}} - 818240\,{u_{5}}\,{u_{8}}\,{u_{4}}^{5}\,{u_{2}}^{2}\,{
u_{3}} + 4322715\,{u_{5}}\,{u_{6}}\,{u_{8}}\,{u_{3}}^{3}\,{u_{4}}
^{3}\,t \\
\mbox{} - 2883690\,{u_{5}}\,{u_{7}}^{2}\,{u_{3}}^{3}\,{u_{4}}^{3}
\,t + 1180\,{u_{7}}\,{u_{8}}^{2}\,{u_{5}}\,{u_{3}}^{4}\,t^{2} - 
1481230\,{u_{5}}^{3}\,{u_{7}}^{2}\,{u_{3}}\,t^{2}\,{u_{4}}^{2}
 \\
\mbox{} - 1237\,{u_{7}}^{3}\,{u_{8}}\,t^{2}\,{u_{2}}^{2}\,{u_{4}}
\,{u_{3}} - 2026175\,{u_{6}}\,{u_{7}}^{2}\,{u_{2}}^{3}\,{u_{4}}^{
2}\,{u_{3}}^{2} + 366130\,{u_{7}}^{2}\,{u_{8}}\,{u_{4}}^{4}\,t^{2
}\,{u_{2}} \\
\mbox{} - 1445980\,{u_{5}}^{2}\,{u_{8}}\,{u_{4}}^{5}\,t\,{u_{2}}
 - 703440\,{u_{6}}^{2}\,{u_{8}}\,{u_{3}}^{4}\,{u_{4}}^{2}\,t + 
436\,{u_{7}}^{4}\,{u_{5}}\,t^{2}\,{u_{2}}^{2}\,{u_{3}} \\
\mbox{} - 566\,{u_{5}}^{2}\,{u_{7}}^{3}\,{u_{6}}\,t^{3}\,{u_{3}}
 - 1759\,{u_{5}}\,{u_{6}}^{2}\,{u_{8}}^{2}\,t^{2}\,{u_{2}}^{2}\,{
u_{3}} + 733440\,{u_{5}}^{5}\,{u_{6}}\,{u_{7}}\,t^{2}\,{u_{2}}
 \\
\end{array}
$$
\newpage 
$$
\begin{array}{l}
\mbox{} - 744856\,{u_{6}}^{3}\,{u_{8}}\,{u_{2}}^{5}\,{u_{4}} - 
354116\,{u_{6}}^{3}\,{u_{8}}\,{u_{2}}^{4}\,{u_{3}}^{2} - 591664\,
{u_{6}}^{3}\,{u_{7}}\,{u_{2}}^{3}\,{u_{3}}^{3} \\
\mbox{} - 725440\,{u_{7}}^{3}\,{u_{6}}\,{u_{2}}^{5}\,{u_{3}} + 
219205\,{u_{5}}^{2}\,{u_{6}}^{2}\,{u_{2}}^{3}\,{u_{4}}^{3} - 
690920\,{u_{5}}^{3}\,{u_{8}}\,{u_{3}}^{5}\,{u_{2}} \\
\mbox{} + 740750\,{u_{7}}^{3}\,{u_{5}}\,{u_{2}}^{4}\,{u_{3}}^{2}
 - 1499205\,{u_{7}}^{3}\,{u_{5}}\,{u_{2}}^{5}\,{u_{4}} + 3120390
\,{u_{5}}^{4}\,{u_{2}}\,{u_{3}}^{2}\,{u_{4}}^{3} \\
\mbox{} - 368395\,{u_{7}}^{4}\,t\,{u_{2}}^{4}\,{u_{4}} - 1299270
\,{u_{5}}^{3}\,{u_{6}}\,{u_{3}}^{5}\,{u_{4}} - 90\,{u_{8}}^{3}\,{
u_{6}}\,{u_{2}}^{5}\,t \\
\mbox{} + 1859340\,{u_{6}}^{4}\,{u_{2}}^{2}\,{u_{4}}^{3}\,t + 
1035412\,{u_{6}}^{4}\,{u_{2}}^{3}\,{u_{3}}^{2}\,{u_{4}} - 287910
\,{u_{6}}^{3}\,{u_{5}}\,{u_{3}}^{5}\,{u_{2}} \\
\mbox{} - 24246\,{u_{6}}^{3}\,{u_{7}}\,{u_{3}}^{5}\,t + 184510\,{
u_{7}}^{4}\,t^{2}\,{u_{2}}^{2}\,{u_{4}}^{2} + 368610\,{u_{7}}^{4}
\,t\,{u_{2}}^{3}\,{u_{3}}^{2} \\
\mbox{} + 11583\,{u_{6}}^{4}\,{u_{3}}^{4}\,t\,{u_{4}} + 787540\,{
u_{6}}^{2}\,{u_{7}}^{2}\,{u_{2}}^{5}\,{u_{4}} + 184050\,{u_{6}}^{
2}\,{u_{8}}^{2}\,{u_{3}}^{4}\,t^{2} \\
\mbox{} - 320\,{u_{6}}^{2}\,{u_{8}}^{2}\,{u_{4}}^{3}\,t^{3} + 
518266\,{u_{5}}^{2}\,{u_{7}}^{2}\,{u_{2}}^{2}\,{u_{3}}^{4} + 
743040\,{u_{5}}^{2}\,{u_{7}}^{2}\,{u_{4}}^{4}\,t^{2} \\
\end{array}
$$
\newpage 
$$
\begin{array}{l}
\mbox{} + 340\,{u_{7}}\,{u_{8}}^{2}\,t^{2}\,{u_{2}}\,{u_{4}}^{3}
\,{u_{3}} + 99795\,{u_{6}}\,{u_{7}}\,{u_{8}}\,{u_{3}}^{3}\,t\,{u
_{4}}^{2}\,{u_{2}} + 1413765\,{u_{6}}\,{u_{7}}\,{u_{8}}\,{u_{3}}
\,t\,{u_{4}}^{3}\,{u_{2}}^{2} \\
\mbox{} - 367254\,{u_{7}}^{4}\,{u_{3}}^{2}\,t^{2}\,{u_{4}}\,{u_{2
}} - 739962\,{u_{6}}^{4}\,{u_{3}}^{2}\,t\,{u_{4}}^{2}\,{u_{2}} + 
7965345\,{u_{5}}\,{u_{6}}^{2}\,{u_{7}}\,{u_{3}}^{2}\,t\,{u_{4}}^{
2}\,{u_{2}} \\
\mbox{} + 175\,{u_{8}}^{3}\,{u_{5}}\,t^{2}\,{u_{2}}^{2}\,{u_{4}}
\,{u_{3}} + 36330\,{u_{6}}^{3}\,{u_{7}}\,{u_{3}}^{3}\,t\,{u_{4}}
\,{u_{2}} - 1536725\,{u_{6}}^{3}\,{u_{5}}\,{u_{3}}\,t\,{u_{4}}^{3
}\,{u_{2}} \\
\mbox{} - 680\,{u_{7}}^{3}\,{u_{8}}\,t\,{u_{2}}^{4}\,{u_{3}} + 
4502365\,{u_{5}}^{2}\,{u_{6}}^{2}\,t\,{u_{2}}\,{u_{4}}^{4} + 
6826864\,{u_{5}}^{2}\,{u_{6}}^{2}\,{u_{2}}^{2}\,{u_{3}}^{2}\,{u_{
4}}^{2} \\
\mbox{} + 1402704\,{u_{5}}^{3}\,{u_{8}}\,{u_{2}}^{3}\,{u_{4}}^{2}
\,{u_{3}} + 4326150\,{u_{5}}^{3}\,{u_{8}}\,{u_{2}}^{2}\,{u_{3}}^{
3}\,{u_{4}} - 196\,{u_{8}}^{3}\,{u_{6}}\,t^{2}\,{u_{2}}^{3}\,{u_{
4}} \\
\mbox{} - 56\,{u_{8}}^{3}\,{u_{6}}\,t^{3}\,{u_{2}}\,{u_{4}}^{2}
 + 94\,{u_{8}}^{3}\,{u_{6}}\,t^{2}\,{u_{2}}^{2}\,{u_{3}}^{2} + 
738420\,{u_{6}}^{3}\,{u_{8}}\,t^{2}\,{u_{2}}\,{u_{4}}^{3} \\
\mbox{} + 747772\,{u_{6}}^{3}\,{u_{8}}\,t\,{u_{2}}^{2}\,{u_{3}}^{
2}\,{u_{4}} - 3735664\,{u_{6}}^{3}\,{u_{7}}\,t\,{u_{2}}^{2}\,{u_{
4}}^{2}\,{u_{3}} - 366789\,{u_{6}}^{3}\,{u_{8}}\,{u_{3}}^{2}\,t^{
2}\,{u_{4}}^{2} \\
\mbox{} + 4380240\,{u_{5}}^{3}\,{u_{6}}\,{u_{3}}\,t\,{u_{4}}^{4}
 + 1402704\,{u_{6}}^{3}\,{u_{5}}\,{u_{3}}^{3}\,t\,{u_{4}}^{2} + 
1503169\,{u_{5}}\,{u_{6}}^{2}\,{u_{7}}\,{u_{2}}^{3}\,{u_{3}}^{2}
\,{u_{4}} \\
\end{array}
$$
\newpage 
$$
\begin{array}{l}
\mbox{} - 2208929\,{u_{5}}\,{u_{6}}^{2}\,{u_{7}}^{2}\,t^{2}\,{u_{
2}}\,{u_{4}}\,{u_{3}} - 755926\,{u_{5}}^{2}\,{u_{6}}^{2}\,{u_{7}}
\,t\,{u_{2}}\,{u_{3}}^{3} \\
\mbox{} - 429928\,{u_{5}}^{2}\,{u_{6}}^{2}\,{u_{7}}\,{u_{2}}^{4}
\,{u_{3}} + 24\,{u_{6}}^{2}\,{u_{7}}^{2}\,{u_{8}}\,{u_{3}}^{2}\,t
^{3} - 729916\,{u_{5}}\,{u_{6}}^{2}\,{u_{7}}\,{u_{8}}\,{u_{2}}^{4
}\,t \\
\mbox{} - 24620\,{u_{6}}^{3}\,{u_{5}}\,{u_{7}}\,t^{2}\,{u_{2}}\,{
u_{4}}^{2} - 728690\,{u_{5}}^{3}\,{u_{8}}\,{u_{7}}\,t^{2}\,{u_{2}
}\,{u_{4}}^{2} - 1475585\,{u_{5}}^{3}\,{u_{8}}\,{u_{7}}\,{u_{3}}
^{2}\,t^{2}\,{u_{4}} \\
\mbox{} - 737800\,{u_{7}}^{3}\,{u_{5}}\,{u_{6}}\,t^{2}\,{u_{2}}^{
2}\,{u_{4}} - 2195695\,{u_{5}}^{3}\,{u_{8}}\,{u_{6}}\,{u_{3}}\,t
^{2}\,{u_{4}}^{2} \\
\mbox{} - 1462540\,{u_{5}}^{2}\,{u_{6}}^{2}\,{u_{8}}\,t^{2}\,{u_{
2}}\,{u_{4}}^{2} - 415\,{u_{6}}^{2}\,{u_{7}}^{2}\,{u_{8}}\,t^{2}
\,{u_{2}}^{3} \\
\mbox{} - 4470273\,{u_{5}}^{2}\,{u_{6}}^{2}\,{u_{7}}\,t\,{u_{2}}
^{2}\,{u_{3}}\,{u_{4}} + 3542894\,{u_{5}}^{3}\,{u_{7}}^{2}\,t\,{u
_{2}}^{2}\,{u_{3}}\,{u_{4}} \\
\mbox{} + 7249510\,{u_{5}}^{3}\,{u_{6}}\,{u_{7}}\,t\,{u_{2}}^{2}
\,{u_{4}}^{2} + 770\,{u_{6}}^{3}\,{u_{7}}\,{u_{8}}\,t^{2}\,{u_{2}
}^{2}\,{u_{3}} - 1496884\,{u_{6}}^{3}\,{u_{5}}\,{u_{8}}\,t\,{u_{2
}}^{3}\,{u_{3}} \\
\mbox{} - 1456521\,{u_{6}}^{3}\,{u_{5}}\,{u_{7}}\,t\,{u_{2}}^{2}
\,{u_{3}}^{2} + 735042\,{u_{7}}^{3}\,{u_{5}}\,{u_{6}}\,t^{2}\,{u
_{2}}\,{u_{3}}^{2} \\
\end{array}
$$
\end{document}